\documentclass[12pt]{amsart}

\usepackage{amsmath,amssymb,amscd,amsthm,amsxtra}

\newtheorem{Thm}{Theorem}[section]

\newtheorem{Cor}[Thm]{Corollary}
\newtheorem{Lem}[Thm]{Lemma}
\newtheorem{Prop}[Thm]{Proposition}
\newtheorem{Def}[Thm]{Definition}

\newcommand{\co}{\text{co\,}}
   
\def\ldots{\mathinner{\ldotp\ldotp\ldotp}}
\def\ldots{\mathinner{\cdotp\cdotp\cdotp}}

\def \cal{\mathcal}
\def \Cal{\mathcal}
\def \Bbb{\mathbb}

\def \spn{\text{span }}
\def\Ave{\mathop{\text {Ave}}}

\begin{document}

\title{Polynomial approximation on convex subsets of
$\mathbb R^n.$
}
\author{Y.A. Brudnyi}
\address {
Department of Mathematics\\
Technion \\
Haifa 32000\\
Israel }
\email{ ybrudnyi@techunix.technion.ac.il}

\author{N. J. Kalton}
\address{Department of Mathematics \\
University of Missouri-Columbia \\
Columbia, MO 65211
}

\email{nigel@math.missouri.edu}
\thanks {The first author was supported by BSF grant 10004; the second
author was supported by NSF grant DMS-9500125}

\begin{abstract}

Let $K$ be a closed bounded convex subset of $\Bbb R^n$; then by a
result of the first author, which extends a classical theorem of Whitney
there is a constant
$w_m(K)$ so that
for every continuous function $f$ on $K$ there is a polynomial $\varphi$
of degree at most $m-1$ so that
$$ |f(x)-\varphi(x)|\le w_m(K)\sup_{x,x+mh\in K}
|\Delta_h^m(f;x)|.$$  The aim of
this paper
is to study the constant $w_m(K)$ in terms of the dimension $n$ and the
geometry of
$K.$  For example we show that $w_2(K)\le \frac12[\log_2n]+\frac54$ and
that for suitable $K$ this bound is almost attained.  We place special
emphasis on the case when $K$ is symmetric and so can be identified as
the unit ball of finite-dimensional Banach space; then there are
connections between the behavior of $w_m(K)$ and the geometry
(particularly the Rademacher type) of the underlying Banach space.  It is
shown for example that if $K$ is an ellipsoid then $w_2(K)$ is bounded,
independent of dimension,
and $w_3(K)\sim \log n.$  We also give estimates for $w_2$ and $w_3$ for
the unit ball of the spaces $\ell_p^n$ where $1\le p\le \infty.$
\end{abstract}

\maketitle

\section{Introduction} \label{intro}
 \setcounter{equation}{0}

{\bf Basic definitions. }
Let $K$ be a closed  subset of $\Bbb R^n$ and  $\Cal P_m$
denote the space of polynomials of total degree at most $m$. If
$f$ is a continuous function on $K$ we set $$E_m(f;K):=\inf_{\varphi\in
\Cal P_{m-1}}\max_{x\in K}|f(x)-\varphi(x)|$$
and
$$ \omega_m(f)=\omega_m(f;K):=\sup_{x,x+h,\ldots,x+mh\in
K}|\Delta_h^m(f;x)|$$ where
$$ \Delta_h^m(f;x):=\sum_{j=0}^m(-1)^{m-j}\binom{m}{j}f(x+jh).$$

We then define the {\it  Whitney constant $w_m(K)$} by:

\begin{equation}\label{e1.1}
w_m(K):=\sup\{E_m(f):\ f\in C(K) \text{ and } \omega_m(f)\le 1\}.
\end{equation}

We will mainly be interested in the case when $K$ belongs to the class
${\cal C}_b(\Bbb R^n)$ of {\it bounded convex subsets} of $\Bbb R^n$  or
to the subclass ${\cal SC}_b(\Bbb R^n)$ of all centrally symmetric
convex subsets of $\Bbb R^n.$  In the latter case $K$ can be identified
with the closed unit ball $B_X$ of an $n$-dimensional Banach space $X$
and it is natural to write $w_m(X)$ in place of $w_m(B_X)$.
As we do not consider unbounded $K$ except in the introduction this
notation does not lead to any ambiguity.

We also define the {\it global Whitney constant} by:
\begin{equation}\label{e1.2}
w_m(n):= \sup\{w_m(K):\ K\in{\cal C}_b (\Bbb R^n)\}.
\end{equation}

In the spirit of the classical paper of Whitney \cite{30} who considers
the case of dimension one\footnote{in this case $w_m(1)=w_m([0,1]).$},
let us consider also the constants $w_m^*(n)$ and $w_m^{**}(n)$ defined
by (\ref{e1.1}) with $K:= \Bbb R^n_{+}:=\{x\in\Bbb R^n:\ x_i\ge 0\}$ and
$K:=\Bbb R^n$ respectively.  Using the techniques of Beurling (cf.
\cite{30}) it is easy to prove the following estimates:
\begin{equation}\label{e1.3}
w^*_m(n)\le 2, \hskip.5truein w_m^{**}(n)\le \min_{1\le j\le m}
1/\binom{m}{j}
\end{equation}

In contrast, the estimates for $w_m(n)$ are not independent of
dimension, and in fact $\lim_{n\to\infty}w_m(n)=\infty$ if $m\ge 2.$

The main goal of this paper is to give ``good" quantitative estimates
for $w_m(n)$ and for $w_m(K)$ in terms of the geometry of the set $K.$
\smallskip

{\bf Remarks.} (a) The inequalities (\ref{e1.3}) are relatively precise.
For
instance $w_2^*(2)\ge 1.$   Concerning the sharpness of the
second inequality even for $n=1$ see \cite{30}.  In fact the Beurling
method yields the more general inequality $w_m(K)\le 2$ provided $K$
satisfies the {\it unbounded cone condition.}  This condition means that
there is an unbounded cone $C$ with vertex at the origin so that
$ K+C\subset K.$

(b) The asymptotic behavior of Whitney's constants does not change if
the supremum in (\ref{e1.2}) is taken over {\it all} convex subsets of
$\Bbb R^n.$  Actually let $\tilde w_m(n):= \sup w_m(K)$ where $K$ runs
over all unbounded convex subsets of $\Bbb R^n.$  Then $w_m(n-1)\le
\tilde w_m(n)$ while compactness arguments show that $\tilde w_m(n)\le
w_m(n).$

(c) If we let
\begin{equation}
\label{1.4}
w_m^{(s)}(n):=\sup_{\dim X=n}w_m(X),
\end{equation}
then $w_m^{(s)}(n)\le w_m(n).$  In the case $m=2$ we have $w_2(n)\le
Cw_2^{(s)}(n)$ for some universal constant $C$ independent of
dimension.  However we do not know of a similar inequality when $m>2.$

(d) In his paper \cite{31} Whitney also proved the finiteness of similar
constants in
 a more general situation in which $C[0,1]$ is replaced by the space
$B[0,1]$ of bounded (not necessarily measurable) functions.  He also
posed the problem for the space $L_0[0,1]$ of measurable functions.
Let us denote by $w_m(K;B)$, (respectively $w_m(K;L_0)$) the
corresponding constants defined by (\ref{e1.1}) allowing $f$ to be
bounded
(respectively, measurable).  One can then prove the inequality:
\begin{equation}\label{e1.5}
w_m(K;B) \le (2^{2m}-1)w_m(K)+2^m \end{equation}
A similar inequality holds for $w_m(K;L_0).$  Since we do not use this
inequality we will omit its proof.

\smallskip

{\bf Prior results: the one-dimensional case.} In \cite{30} Whitney
proved that $w_m(1)<\infty$ for all $m$ and gave numerical estimates for
$w_m(1)$ when $m\le 5.$ Using a different approach, the first-named
author
proved the analog of the Whitney inequality for translation-invariant
Banach lattices and gave, in particular, an effective but rather rough
estimate of $w_m(1)$ for all $m.$  This estimate was subsequently
improved by a research team (K. Ivanov, Binev and Takev) headed by Sendov
who finally showed that $w_m(1)\le 6$ for all $m.$  The most recent
result is due to Kryakin who proved that $w_m(1)\le 2$ for all $m$ (see
\cite{16} for the references).  The only known precise result is
$w_2(1)=\frac12.$
\smallskip

{\bf Prior results: the multidimensional case.} In 1970, the first-named
author \cite{2} established the multidimensional analog of Whitney's
result for translation-invariant Banach lattices.  From this it follows,
in particular, that $w_m(n)<\infty$ for every $m,n.$  Later in a lecture
at Moscow State University he established an estimate $w_2(n)\le C\log
(n+1).$  Following this lecture S. Konyagin suggested that
$w_2(X)<\infty$ for every {\it infinite-dimensional} Banach space $X$
belonging to the class $\cal K$ introduced by the second author
(cf. \cite{9}).  In particular this implies that $w_2(\ell_p^n)$ is
bounded by a constant independent of dimension if $1<p\le \infty.$  This
important observation led to the authors' collaboration on the current
paper.
\smallskip

{\bf Discussion of the main results.} Our main results concern the
Whitney
constants for $m=2$ and $m=3$ (see Section 5 for some results when
$m>3$).  In Section 3, we give a fairly precise estimate for $w_2(n)$,
i.e.,
$$ \frac12\log_2(\left[\frac{n}{2}\right]+1)\le w_2(n)\le
\frac12[\log_2n]+\frac54.$$
Curiously enough $w_2(n)$ is almost attained not for the unit simplex
$S^n$ but for its Cartesian square.  Meanwhile for $S^n$ we prove in
Theorem \ref{3.7} the precise asymptotics are given by
$$ \lim_{n\to\infty}\frac{w_2(S^n)}{\log_2n}=\frac14.$$

We also consider in this section the problem of estimating
$w_2(\ell_p^n)$ for $1\le p\le \infty.$  In particular, we show in
Theorem \ref{3.8} that  $w_2(\ell_1^n)\sim\log n$ while
$\gamma(p):=\sup_n w_2(\ell_p^n)$ is finite for $1<p\le \infty.$  More
precisely $\gamma(p)$ is equivalent up to a logarithmic factor to
$(p-1)^{-1}$ when $p\downarrow 1$; surprisingly, for $2\le p\le \infty,$
the constant
$\gamma(p)$ is bounded by an absolute constant.  This striking difference
in asymptotic behavior is explained by Theorem \ref{3.11} which gives an
upper estimate of $w_2(X)$ in terms of the type $p$ constant $T_p(X)$ of
$X.$

In Section 4, we consider the problem of quadratic approximation on
symmetric convex bodies.  In particular we show in Theorem \ref{4.1}
that
$$ c_1\sqrt{n} \le w_3^{(s)}(n)\le c_2\sqrt{n}\log(n+1)$$ for some
absolute constants $0<c_1,c_2<\infty.$  As in the linear case, however,
better estimates are available for $\ell_p^n-$spaces.  Our results in
this case are consequences of Theorem \ref{4.5} giving upper and lower
estimates of $w_3(X)$ by the type 2 constant of $X$, $T_2(X)$ and the
cotype
2 constant of $X^*$, $C_2(X^*).$  Actually we show that
$$c_1C_2(X^*)^{-8} \le \frac{w_3(X)}{\log(n+1)}\le c_2 T_2(X)^2$$ for
some absolute positive constants $c_1,c_2.$    As corollaries of this
inequality we have, for example, that $w_3(\ell_2^n)\sim\log (n+1)$
while
$$ c_1\log(n+1)\le w_3(\ell_{\infty}^n)\le c_2(\log(n+1))^2$$
with absolute constants $c_1,c_2>0$.  See Theorem \ref{4.3}.

In Section 5 we discuss a few  estimates for $w_m^{(s)}(n)$
and
$w_m(\ell_p^n)$ with $m\ge 4$.  In particular, we prove in Corollary 5.5,
that
\begin{equation}
\label{1.6}
w_m^{(s)}(X)\le cn^{\frac{m}{2}-1}\log(n+1)
\end{equation}
where $n=\dim X$ and $C$ is an absolute constant.  Once again for
$\ell_p^n$-spaces we have better estimates. For instance,
$$ w_m(\ell_p^n)\le Cn^{\frac{m-3}{2}}\log(n+1)$$
for $2\le p\le \infty$ and $m\ge 3$ while $w_m(\ell_1^{n})\sim\log(n+1).$
Particularly striking is the fact that there is a dimension-free upper
bound for $w_m(\ell_p^n)$ for (fixed) arbitrary $m$ if $0<p<1$ as in
Theorem \ref{5.8}.

Our arguments depend, in part, on some deep results of the local theory
of Banach spaces.  Most of them are concentrated in the proofs of
Theorems \ref{3.7} and \ref{4.5}.  We also need a refinement of the main
result Theorem 1.1 of the paper \cite{11} and a version of Maurey's
extension principle \cite{20} using a dual cotype 2 assumption
in place of
the usual type 2 assumption.  The proof of the first result is
presented
in Section 3 while the required ingredients of the proof are presented in
Section 2.  This section also contains the proof of the second result and
those of two results related to the {\it homogeneous} versions of the
Whitney's constants.

Let us discuss our results in connection with the {\it curse of
dimension}, which, roughly speaking asserts that the
computational complexity of a function of $n$ variables grows
exponentially in $n.$  In situations where this can be precisely
formulated and proved it is, in general, a statement of the complexity of
a universal (e.g. linear) approximation method for functions in a  given
class.  It may be anticipated that approximation methods for individual
functions can be much more efficient.  In these terms we can consider
$w_m(K)$ as a measure of approximation of $f\in C(K)$
satisfying $\omega_m(f;K)\le 1$ by polynomials of degree $m-1.$  We can
then compare $w_m(K)$ with a linearized Whitney constant $w_m^l(K)$ which
is defined by
$$ w_m^l(K)= \inf_{L}\sup_{\omega_m(f:K)\le 1}\|f-Lf\|_{C(K)}$$
where $L$ runs through all linear operators $L:C(K)\to\Cal P_{m-1}.$
In the case when $K=B_{\ell_2^n}$ this quantity has been estimated by
Tsarkov \cite {29}, who proves
$$ w_m^l(K) \sim n^{(m-1)/2}.$$
Our results show that $w_m(K)\le Cn^{(m-3)/2}\log(n+1)$ for $m\ge 3.$
Thus we have a marked improvement over linear methods which is especially
striking when $m=3$ since $w_3^l(K)\sim n$ but $w_3(K)\sim \log (n+1).$
\smallskip

{\bf Remarks on the infinite-dimensional case.} There is an obvious
generalization of the Whitney constant $w_m(X)$ to the case when $X$ is
an infinite-dimensional Banach space (or even quasi-Banach space).  In
this case it is quite possible that $w_m(X)=\infty.$  Let us consider
first the case when $m=2.$  We recall (cf. \cite{9} or \cite{13}) that a
Banach space $X$ is called a $\cal K$-space if, whenever $f:X\to \Bbb R$
is a quasilinear map (see Section 2) then there is a linear functional
$g:X\to\Bbb R$ with $\sup\{|f(x)-g(x)|:\ x\in
B_X\}:=\|f-g\|_{B_X}<\infty.$ There is a clear connection
between the above condition and $w_2(X)<\infty.$  However, since the
definition of a $\cal K$-space allows for discontinuous $f$ (and $g$) it
is not clear that these conditions are equivalent.   They are equivalent
if $X$ has the bounded approximation property.

For the case $m\ge 3$ it is possible to show that $w_m(X)=\infty$ for
most classical spaces..  More precisely, $w_m(X)=\infty$ for $m\ge 3$ if
$X$ contains uniformly complemented $\ell_p^n$'s for some $1\le p\le
\infty$; this includes the case when $X$ has nontrivial type.  The same
conclusion can also be reached if $X^*$ has cotype 2 and this covers the
case of the space $P$ constructed by Pisier \cite{24} as an example
of
a space which does not contain uniformly complemented finite-dimensional
subspaces.

For infinite-dimensional quasi-Banach spaces, this situation is quite
different.  For example $w_m(\ell_p)<\infty$ for any $m\in\Bbb N$ and
$0<p<1.$  The case of $L_p(0,1)$ is even more remarkable, since
 $w_m(L_p)<\infty$ for every $m\in\Bbb N$  and yet the only polynomials
on $L_p$ are constant (because $L_p$ has trivial dual).  Thus if
$F:B_{L_p}\to\Bbb R$ is continuous, satisfies $F(0)=0$ and
$\omega_m(F)\le
1$ then $\|F\|_{B_{L_p}}\le C$ where $C=C(m,p).$

It is worth perhaps remarking that although the paper does not explicitly
use the theory of twisted sums of Banach and quasi-Banach spaces, this
theory is implicit in many of the results, and there is a clear
connection with ideas in \cite{9}, \cite{KP}, \cite{12} and \cite{Rib}.

\smallskip

{\bf The stability of the equation $\Delta^m_hf=0.$}  There is an
alternative viewpoint for the results presented in this paper.  It is
well-known that a continuous function $f$ defined on a convex set $K$ is
a polynomial of degree $m-1$ if and only if $f$ satisfies the functional
equation $\Delta^m_hf=0.$  So the Whitney constant $w_m(K)$ can be
regarded as a measure of stability of this equation.  Stability problems
of this type go back to the work of Hyers and Ulam.  We note in this
connection the work of Casini and Papini \cite{3} and a recent preprint
of Dilworth, Howard and Roberts \cite{4} on stability of convexity
conditions.

{\bf Conjectures.}  The work in this paper was motivated by certain
conjectures, and it may be helpful to list them here.

1.)  If $m\ge 2$ then
$$ w_m(n) \sim w_m^{(s)}(n)\sim n^{\frac{m}{2}-1}\log (n+1)$$ as $n\to
\infty.$

This conjecture is proved for $m=2$ while the upper estimate for
$w_m^{(s)}(n)$ is established for all $m\ge 2.$  As the lower bound for
$m\ge 3$ we have only the inequalities $w_m(n)\ge w_m^{(s)}(n)\ge c\sqrt
n.$

2.) If $m\ge 3$ and $1\le p<\infty$ then
$$ w_m(\ell_p^n)\sim\log (n+1)$$
as $n\to\infty.$

This result is established for $p=1$ and for $m=3$ and $2\le p<\infty$
while the lower bound is established for all $m\ge 3.$  It is quite
possible that this conjecture is way off the mark when $m\ge 4.$

3.) $w_2(\ell_{\infty}^n)$ is ``small.''  We propose the conjecture that
$w_2(\ell_{\infty}^n)\le 2$ for all $n.$  The only known results are
$w_2(\ell_{\infty}^1)=\frac12$ and $w_2(\ell_{\infty}^2)=1.$  Note that
if our conjecture were to hold then for every convex function $f$ on the
$n$-cube $Q^n$ we would have the inequality $E_2(f;Q^n)\le
\omega_2(f:Q^n).$
\bigskip

4.) If $X$ is an infinite-dimensional Banach space then $w_3(X)=\infty.$

\section{Preliminary results}\label{prelim}
\setcounter{equation}0

{\bf Homogeneous Whitney constants.} Suppose that $X$ is an
$n$-dimensional Banach space.  We consider the homogeneous version of the
Whitney problem.  We say that a function $f:X\to\Bbb R$ is
{\it $m$-homogeneous} if $f(ax)=a^mf(x)$ whenever $a\in\Bbb R$ and $x\in
X.$

\begin{Def}\label{2.1} The homogeneous Whitney constant $v_m(X)$ for
$m\ge 2$ is
the least constant so that if $f$ is an $(m-1)$-homogeneous  continuous
function on $X$ there is an $(m-1)$-homogeneous polynomial $\varphi$ so
that for all $x\in X$,

\begin{equation}\label{e2.1}
|f(x)-\varphi(x)|\le v_m(X)\|x\|^{m-1}\omega_m(f),\end{equation}
where $\omega_m(f)=\omega_m(f;B_X).$
 \end{Def}

If $f$ is continuous and homogeneous (i.e. 1-homogeneous) then
$$ |f(x+y)-f(x)-f(y)|\le \omega_2(f;B_X)\max(\|x\|,\|y\|).$$
Thus $f$ is {\it quasilinear} in the sense of \cite{9}.  This connection
was first noticed by S. Konyagin and the following result is essentially
due to him (see remarks in the introduction):

\begin{Prop}\label{2.2} If $X$ is a finite-dimensional normed space then
$$ v_2(X)\le w_2(X)\le 4v_2(X)+\frac32.$$ \end{Prop}

\begin{proof}
 If $f:X\to \Bbb R$ is continuous and homogeneous,
then an affine function of best approximation on the ball can be taken as
a linear functional, $x^*$ say, and then $|f(x)-x^*(x)|\le
w_2(X)\omega_2(f;B_X)$ so that $v_2(X)\le w_2(X).$

Conversely, suppose $f:B_X\to \Bbb R$ is continuous and that
$\omega_2(f)\le 1.$
Let us note  that any $x,y\in B_X$ we have
\begin{equation}\label{line}
|f(tx+(1-t)y)-tf(x)-(1-t)f(y)|\le 2E_1(f;[x,y])\le 1.\end{equation}
This follows from applying Whitney's one-dimensional result to the
line-segment $[x,y]$, since $w_2(1)=\frac12.$

We define
$g$  on
$X$ by  $g(x)=\frac12\|x\|(f(x/\|x\|)-f(-x/\|x\|))$ for $x\neq 0$ and
$g(x)=0.$  Then $g$ is continuous and homogeneous.  We will show first
that
$\omega_2(g;B_X)\le 4.$

Suppose $x,y\in X$ are not both zero.  Let
$$ \lambda=\frac{\|x\|}{\|x\|+\|y\|},\ \mu=\frac{\|y\|}{\|x\|+\|y\|},\
\nu=\frac{\|x+y\|}{\|x\|+\|y\|}$$ and choose $u,v,w\in B_X$ so that
$\|u\|=\|v\|=\|w\|=1$ and
$$ \|x\|u=x, \qquad \|y\|v=y, \qquad \text{and } \|x+y\|w=x+y.$$
Then for $\epsilon=\pm1$,
$$ I_{\epsilon}:=|f(\epsilon(\lambda u+\mu v))-\lambda f(\epsilon
u)-f(\epsilon v)|\le 1$$
by applying (\ref{line}). Similarly
$$ J_{\epsilon}:= |f(\epsilon(\lambda u+\mu v))-\nu f(\epsilon w)-(1-\nu)
f(0)|\le 1.$$
From the definition of $g$ we have:
\begin{equation}\begin{align*}
|g(x)-2g(\frac12(x+y))+g(y)|&=|g(x)-g(x+y)+g(y)|\\
&\le \frac12\|x+y\|\sum_{\epsilon=\pm 1}(I_{\epsilon}+J_{\epsilon})\\
&\le 2\|x+y\|\le 4.
\end{align*}\end{equation}
Hence $\omega_2(g:B_X)\le 4.$

This implies that there exists $x^*\in X^*$ so that if $\|x\|\le 1,$
$$ |g(x)-x^*(x)|\le 4v_2(X).$$   We will choose $\varphi(x)=x^*(x)+f(0)$
as an affine approximation to $f.$  If $\|x\|=1$ then,
\begin{equation}\begin{align*}
 |f(x)-\varphi(x)|&\le 4v_2(X) +|f(x)-f(0)-g(x)| \\
&\le 4v_2(X)+\frac12 |f(x)+f(-x)-2f(0)|\le
4v_2(X)+\frac12.\end{align*}\end{equation}

Now suppose $\|y\|\le 1.$  We write $y=tx$ where $\|x\|=1$ and $0\le t\le
1.$  By (\ref{line}) we have:
$$ |f(y)-tf(x)-(1-t)f(0)|\le 1$$
and hence
$$ |f(y)-\varphi(y)| \le 4v_2(X)+\frac{3}{2}.$$
This completes the proof.\end{proof}

The following Lemma gives a uniform estimate on $w_m(X)$ for all $X$ of
dimension $n$ (cf. \cite{2}):

\begin{Lem}\label{2.3}
For any $m\ge 2$, and any $n$-dimensional Banach space $X,$
$$w_m(X) \le 2+T_{m-1}(\sqrt n)(2+w_m(\ell_2^n)),$$
where $T_k(t):=\cos(k\arccos t)$ is the Chebyshev polynomial of  degree
$k.$\end{Lem}

\begin{proof}
By a well-known result of John \cite{8} there is a Euclidean norm
$\|\cdot\|_E$ on $X$ so that
$$ n^{-1/2}\|x\|_E \le \|x\|_X\le \|x\|_E$$
for $x\in X.$  Now suppose that $f:B_X\to\Bbb R$ is continuous and
$\omega_m(f)\le 1.$  Restricting $f$ to $B_E$ we can find a polynomial
$\varphi\in \cal P_{m-1}$ with $|f(x)-\varphi(x)|\le w_m(\ell_2^n)$ for
$x\in B_E.$  Fix any $x\in B_X.$  By the definition of the Whitney
constant and Kryakin's theorem \cite{16} there is a polynomial
$\psi\in\cal P_{m-1}(\Bbb R)$ so that
$$|f(tx)-\psi(t)|\le w_m([0,1])\le 2$$ for $|t|\le 1.$  Hence for $|t|\le
n^{-1/2}$ we have
$$ |\varphi(tx)-\psi(t)|\le 2+w_m(\ell_2^n).$$
According to the Chebyshev inequality (see e.g. \cite{Riv} p. 108) it
follows that for $|t|\le 1$
$$ |\varphi(tx)-\psi(t)|\le T_{m-1}(\sqrt n)(2+w_m(\ell_2^n)).$$
The result now follows easily.\end{proof}

Let us also note at this point that essentially the same argument gives
us the following elementary estimate:

\begin{Lem}\label{2.3a} Let $X,Y$ be two $n$-dimensional normed spaces
and let $d:=d(X,Y)$ be the  Banach-Mazur distance between them. Then
$$ w_m(Y)\le 2+T_{m-1}(d)(2+w_m(X))$$ and
$$ v_m(Y)\le d^{m-1}w_m(X).$$
\end{Lem}

\begin{proof} We may suppose that $\|\,\|_Y$ and $\|\,\|_X$ are two norms
on $\Bbb R^n$ so that $d^{-1}\|x\|_X \le \|x\|_Y\le \|x\|_X$ for
$x\in\Bbb R^n.$  The first estimate is proved just as in Lemma \ref{2.3}.
The second estimate follows easily from the definition of $v_m(X)$ using
(\ref{e2.1}).\end{proof}

We now prove a much more general version of Proposition \ref{2.2}.

\begin{Prop}\label{2.4}
Suppose that $m\ge 2.$  Then there is a constant $C=C(m)$ (independent of
$X$)  so that for every finite-dimensional Banach space $X$,
$$ C^{-1}\max_{2\le k\le m}v_k(X)\le w_m(X)\le C\max_{2\le k\le
m}v_k(X).$$
\end{Prop}

\begin{proof} First choose for each $0\le i\le m-1$ real numbers
$(c_{ij})_{j=1}^m$ so that for any polynomial $\varphi$ in one variable
of degree at most $m-1$ we have:
\begin{equation}\label{e2.2}
\frac{\varphi^{(i)}(0)}{i!} =\sum_{j=1}^mc_{ij}\varphi(\frac{j}{m}).
\end{equation}
In particular we have
\begin{equation}\label{e2.3}
\sum_{j=1}^mc_{ij}(\frac{j}{m})^k=\delta_{ik}.
\end{equation}
for $0\le i,k\le m.$  Hence, if $\varphi\in\cal P_{m-1}$ then
$\psi(x):= \sum_{j=1}^mc_{k-1,j}\varphi(\frac{jx}{m})$ is a
$(k-1)-$homogeneous polynomial.

Using this, let us first prove that
\begin{equation}\label{e2.4}
v_k(X)\le C(m)w_m(X), \hskip.5truein 2\le k\le m
\end{equation}
In fact if $f:X\to \Bbb R$ is continuous and $(k-1)-$homogeneous with
$\omega_k(f)=\omega_k(f;B_X)\le 1$ then $\omega_m(f)\le 2^{m-k}$ and so
there exists a polynomial $\varphi\in\cal P_{m-1}$ with
$$ |f(x)-\varphi(x)|\le 2^{m-k}w_m(X).$$
Now $f(x)=\sum_{j=1}^mc_{k-1,j}f(\frac{jx}{m})$
by the $(k-1)-$homogeneity of $f$ and (\ref{e2.3}), and this inequality
leads to the estimate
$$ |f(x)-\psi(x)|\le 2^{m-k}(\sum_{j=1}^m|c_{k-1,j}|)w_m(X)$$
for $x\in B_X$ where $\psi(x):=
\sum_{j=1}^mc_{k-1,j}\varphi(\frac{jx}{m})$ is a $(k-1)$-homogeneous
polynomial.  Hence (\ref{e2.4}) follows.

Conversely let $V:=\max_{2\le k\le m}v_k(X).$  Suppose $f\in C(B_X)$
with
$\omega_m(f)\le 1$.  Then for each $x$ with $\|x\|=1$ and $1\le k\le
m-1$ we define
$g_k(x)=\sum_{j=1}^mc_{kj}f(\frac{jx}{m})$ and extend $g_k$ to be
$k$-homogeneous. It is easy to see that each $g_k$ is continuous.  We
also let $g_0(x)=f(0)$ for all $x\in X.$

By the one-dimensional result \cite{16} for each $x$ with $\|x\|=1$ there
is a polynomial $\varphi$ on $[0,1]$ of degree at most $m-1$ so that
$$ |f(tx)-\varphi(t)|\le 4$$ for $0\le t\le 1.$  Hence
$$ |g_k(x)-\frac{\varphi^{(k)}(0)}{k!}| \le 4\max(1,\sup_{1\le l\le
m-1}\sum_{j=1}^m|c_{lj}|)\le C_1$$
where $C_1=C_1(m).$
Then, for any $x\in B_X$ we have
$$ |f(x)-\sum_{k=0}^{m-1}g_k(x)|\le 4+mC_1=C_2.$$
Using (\ref{e2.3}) for $1\le k\le m-1,$ we have the identity
$$ \sum_{j=1}^mc_{kj}f(\frac{jx}{m})-g_k(x)=\sum_{j=1}^mc_{kj}[f(\frac
{jx}{m})-\sum_{s=0}^{m-1}g_s(\frac{jx}{m})] $$
and we can deduce
$$ |g_k(x)-\sum_{j=1}^mc_{kj}f(\frac{jx}{m})|\le
C_2\sum_{j=1}^m|c_{kj}|\le C_3(m).$$
Hence
\begin{equation}\label{e2.5} \omega_m(g_k)\le
2^mC_3+\sum_{j=1}^m|c_{kj}|=C_4\end{equation} for
$1\le k\le m-1.$

We now deduce from (\ref{e2.5}) that
\begin{equation}\label{e2.6}
\omega_{k+1}(g_k)\le C_5(m)\end{equation} for $1\le k\le m-1.$  Indeed
let $x,x+(k+1)h\in B_X$ and let $F:=\spn \{x,h\}$ be the linear space
generated by $x,h$.  By Lemma \ref{2.3} and the multivariate Whitney type
inequality (in dimension 2) \cite{2} we can find a polynomial $\psi_F$ of
degree at most $m-1$ so that
$$ |g_k(y)-\psi_F(y)|\le C_6\omega_m(g_k)$$
for $y\in B_F$ where $C_6=C_6(m).$  But, arguing as before, we can
replace
$\psi_F$ by $\sum_{j=1}^mc_{kj}\psi_F(\frac{jx}{m})$ and this allows us
to assume that $\psi_F$ is homogeneous of degree $k$ (by similar
arguments to those used above.)  Hence
$$ |\Delta_h^{k+1}g_k(x)|= |\Delta_h^{k+1}(g_k-\psi_F)(x)|\le
2^{k+1}C_6\omega_m(g_k)$$  Combining with (\ref {e2.5}) we get
(\ref{e2.6}).  Then we can conclude that there is a $k$-homogeneous
polynomial $\psi_k$ on $X$ so that
$$ |g_k(x)-\psi_k(x)|\le C_7(m)V$$
for $x\in B_X.$  Finally if we set
$\psi(x)=g_0(x)+\sum_{k=1}^{m-1}\psi_k(x)$ then
$$ |f(x)-\psi(x)| \le (C_2+mC_7) \le C_8(m)V$$ for $\|x\|\le 1$ and so
$w_m(X)\le CV$ for a constant $C$ depending only on $m.$ \end{proof}

\begin{Cor}\label{2.5} If $2\le l\le m$ there is a constant $C=C(l,m)$
so that $$w_l(X)\le C(l,m)w_m(X).$$\end{Cor}

{\bf Remark.} All the above results are clearly true
(with
constants also depending on $r$) for $r$-normed finite-dimensional
spaces.
Recall (cf. \cite{13}) that $\|\cdot\|$ is an $r$-norm on $X$ if we
have
\begin{itemize}
\item{(1)} $\|x\|\ge 0$ with equality if and only if $x=0;$
\item{(2)} $\|ax\|=|a|\|x\|$  for $a\in \Bbb R$ and $x\in X$;
\item{(3)} $\|x_1+x_2\|^r\le \|x_1\|^r+\|x_2\|^r$ for $x_1,x_2\in X.$
\end{itemize}

We note only that in the proof of Lemma \ref{2.3}, John's theorem is
replaced by its $r$-normed generalization due to Peck \cite{22}.
\smallskip

 {\bf Indicators of finite-dimensional Banach lattices}
Let $X=\{{\Bbb R}^{n+1},\|\cdot\|_X\}$ be an $(n+1)-$dimensional Banach
lattice.  In our setting this simply implies that if
$x=(x_i)_{i=1}^{n+1}$ and $y=(y_i)_{i=1}^{n+1}$ with $|x|\le |y|$ (i.e.
$|x_i|\le |y_i|$ for $i=1,2,\ldots,n+1),$ then $\|x\|_X\le \|y\|_X.$

\begin{Def}\label{2.7} (\cite{11}) The {\it indicator} $\Phi_X$ of $X$
is
the function defined on the simplex $S^n:=\{u\in{\Bbb R}^{n+1}:\ u\ge 0,\
\sum_{i=1}^{n+1}u_i=1\}$ by
\begin{equation}\label{e2.7}\Phi_X(u):=\sup_{\|x\|_X\le
1}\sum_{i=1}^{n+1}u_i\log_2|x_i|
\end{equation}
\end{Def}

Here we set $0\log_20=0.$  We remark first that we use logarithms
base two in place of natural logarithms as in \cite{11} for convenience.
We also remark that in \cite{OS}
the same function is called the {\it entropy function} of $X.$

We denote by $\Lambda$ the functional
$\Lambda(u)=\sum_{i=1}^{n+1}u_i\log_2|u_i|$.  Let us note the following
straightforward properties of $\Phi_X.$

\begin{Prop}\label{2.11}(a) $\Phi_{\ell_1^{n+1}}=\Lambda,$\newline
(b) If $a_i>0$ for $1\le i\le n+1$, $1\le p<\infty$  and
$\ell_p^{n+1}(a)$ is defined by the norm
$$ \|x\|_{\ell_p^{n+1}(a)}:=(\sum_{i=1}^{n+1}a_i^p|x_i|^p)^{1/p}$$
then
$$
\Phi_{\ell_p^{n+1}(a)}(u)=\frac1p(\Lambda(u)-\sum_{i=1}^{n+1}u_i\log_2
a_i)$$
(c) If $\|\cdot\|_X$ and $\|\cdot\|_Y$ are $C$-equivalent i.e.
$C^{-1}\|x\|_X\le \|x\|_Y\le C\|x\|_X$ for all $x\in\Bbb R^{n+1}$ then
$$ |\Phi_X(u)-\Phi_Y(u)|\le \log_2C$$
for $u\in S^n.$
\end{Prop}

Let us use $\langle x,y\rangle$ to denote the standard inner-product on
$\Bbb R^{n+1}.$  Then if $X$ is a Banach lattice we define the dual
space $X^*$ by
$$ \|x^*\|_{X^*}:=\sup\{|\langle x^*,x\rangle|:\ \|x\|_X\le 1\}.$$

If $X_0,X_1$ are two $(n+1)-$dimensional Banach lattices we define
the (Calder\'on)  interpolation space
$X_{\theta}=X_0^{1-\theta}X_1^{\theta}$ for $0<\theta<1$ by
$$
\|x\|_{X_{\theta}}:=\inf\{\|x_0\|_{X_0}^{1-\theta}
\|x_1\|_{X_1}^{\theta}\}$$ where the infimum is taken over all
$x_0,x_1\in \Bbb R^{n+1}$ satisfying
$$ |x|\le |x_0|^{1-\theta}|x_1|^{\theta}.$$

The following results are taken from \cite{11}:

\begin{Thm}\label{2.10} (a) For any Banach lattice $X$ on $\Bbb R^{n+1},$
$$\Phi_X+\Phi_{X^*}=\Lambda.$$
(b) If $X_0,X_1$ are two Banach lattices on $\Bbb R^{n+1},$ then
$$
\Phi_{X_0^{1-\theta}X_1^{\theta}}=(1-\theta)\Phi_{X_0}+\theta\Phi_{X_1}.$$
\end{Thm}

Note that (b) is a simple consequence of the definitions, while (a)
follows from the deep duality theorem of Lozanovskii \cite{19} (which is
essentially equivalent to the statement that
$X^{1/2}(X^*)^{1/2}=\ell_2^{n+1}$ for any Banach lattice $X.$
It is not hard to see that $\Phi_X$ is a convex function satisfying
$\delta_2(\Phi_X)\le 1$ where $\delta_2: C(S^{n})\to \Bbb R$ is defined
by
$$ \delta_2(f):=\sup\{|f(\alpha u+(1-\alpha)v)-\alpha
f(u)-(1-\alpha)f(v)|\}$$ where the supremum is taken over all $0\le
\alpha\le 1$ and $u,v\in S^n.$

The main result of \cite{11} gives, in our setting, a form of converse
to this statement.

\begin{Thm}\label{2.8}  For each $0<\epsilon<\frac12$  there is a
constant $C=C(\epsilon)$ so that
 whenever $n\in\Bbb N,$ and $f\in C(S^n)$
satisfies
$\delta_2(f)\le 1-\epsilon$ there is a Banach lattice $X$ so that
$$ |f(u)-(\Phi_X(u)-\Phi_{X^*}(u))|\le C $$
for all $u\in S^n.$
\end{Thm}

One of our goals is to refine this result to give a very general
representation
for functions on $S^n$ in terms of the parameter $\omega_2(f).$  This
will be achieved in Theorem \ref{2.9} below.
\smallskip

{\bf Extension theorems of Maurey type.}
We recall that if $X$ is a Banach space and $1<p\le 2$ then $X$ is
said to have {\it type $p$} if there is a constant $C$ so that for any
$x_1,\ldots,x_n\in X$ we have
$$
\left(\Ave_{\epsilon_i=\pm1}\|\sum_{i=1}^n\epsilon_ix_i\|^p\right)^{1/p}
\le C\left(\sum_{i=1}^n\|x_i\|^p\right)^{1/p}.$$
The best constant $C$ is called the type $p$ constant of $X$ and denoted
by
$T_p(X).$

$X$ is said to have {\it cotype $q$} where $2\le q<\infty$ if there is a
constant $C$ so that
 for any
$x_1,\ldots,x_n\in X$ we have
$$\left(\sum_{i=1}^n\|x_i\|^q\right)^{1/q}\le
C\left(\Ave_{\epsilon_i=\pm1}\|\sum_{i=1}^n\epsilon_ix_i\|^q\right)^{1/q}.
$$ The best such constant is denoted by $C_q(X).$

We remark that if $\dim X=n$ then we have $T_p(X)\le n^{1/p}$ and
$C_q(X)\le n^{1/p}$ where $\frac1p+\frac1q=1.$
We also have a duality relationship, namely $C_q(X^*)\le T_p(X).$

Let $X$ and $Y$ be
finite-dimensional Banach spaces and suppose $E$ is a linear subspace of
$X.$

\begin{Def}\label{2.12} The {\it extension constant} ${\cal E}_X(E,Y)$ is
infimum of all constants $M$ so that every  linear map $T:E\to Y$
has a linear extension $T_1:X\to Y$ with $\|T_1\|\le M\|T\|.$
\end{Def}

The Maurey extension principle \cite{20} gives the following estimate for
$Y=\ell_2^m$:
\begin{equation}
\label{e2.8}
{\cal E}_X(E,\ell_2^{m}) \le T_2(X)
\end{equation}

In order to extend this principle to non-Hilbertian $Y$ we can use
the abstract Grothendieck theorem of Pisier.  This states (\cite{25}
Theorem
4.1) that if $T:E\to Y$ is a linear map then there is a factorization
$V:X\to\ell_2^m$ and $U:\ell_2^m\to Y$ so that $\|U\|\|V\|\le
(2C_2(E^*)C_2(Y))^{3/2}.$
(In fact we can do a little better, i.e. we can obtain $\|U\|\|V\|\le
 CC_2(X^*)C_2(Y)(1+\log C_2(X^*)C_2(Y))$
 Putting these estimates together we obtain
\begin{equation}\label{e2.9}
{\cal E}_X(E;Y)\le (2C_2(E^*)C_2(Y))^{\frac32}T_2(X).
\end{equation}

We will need an analogous result with a cotype assumption on $X^*$ in
place of the type restriction on $X.$  The following result may be known
to specialists but we have not been able to find it in the literature:

\begin{Thm}\label{2.13} There is an increasing function
$\psi:(1,\infty)\to(1,\infty)$ so that
\begin{equation}\label{e2.10}
{\cal E}_X(E,\ell_2^m) \le \psi(T_2(X/E))C_2(X^*).\end{equation}
\end{Thm}

\begin{proof} Suppose that $T_0:E\to\ell_2^m$ with $\|T_0\|\le 1.$ We
need to find an extension $T:X\to \ell_2^m$ of $T_0$ with norm majorized
by the right-hand side of (\ref{e2.10}). To do
this we follow an extension technique of Kisliakov which is used heavily
in \cite{12}.  Consider the space $Z=X\oplus_1\ell_2^m$ i.e. $Z=X\times
\ell_2^m$ algebraically with norm $\|(x,y)\|_Z
=\|x\|_X+\|y\|_{\ell_2^m}.$ Then
$Z^*=X^*\oplus_{\infty}\ell_2^m$ i.e. $Z^*=X^*\times \ell_2^m$ with norm
$\|( x^*,y^*)\|_{Z^*}=\max(\|x^*\|_{X^*},\|y^*\|_{\ell_2^m}).$  Since
$C_2(\ell_2^m)=1$ we have
\begin{equation}\label{e2.11}
C_2(Z^*)\le \sqrt2 C_2(X^*).
\end{equation}

Let $G:=\{(x,-T_0x):\ x\in E\}\subset Z.$  Let $ Y:=Z/G$ and let
$Q:Z\to  Y$ be the quotient map. Note that $Q$ maps $\{0\}\times
\ell_2^m$ isometrically onto a subspace $H$ of  $Y$ and that
by \ref{e2.8} there is a projection $P:Y\to H$ with $\|P\|\le T_2(Y).$
Let $S:X\to Z$ be defined by $S(x):=(x,0).$  Then $PQS$ can be regarded
as an extension of $T_0$; more precisely, $T:=\text{Pr}_2(Q^{-1}PQS)$
extends $T_0$ where $\text{Pr}_2(x,y):=y$ and $Q^{-1}$ is the inverse of
$Q$ on
$\{0\}\times\ell_2^m.$  Then $\|T\|\le \|P\|\le T_2(Y).$  It therefore
remains only to estimate $T_2(Y).$

Fix $1<p<2.$  Note that $Y/H$ is isometric to $X/E$.  Hence by arguments
that go back to the paper \cite{5} (see \cite{10} for details) we have
the estimate $T_p(Y)\le \varphi(T_2(X/E))$ for a suitable increasing
function $\varphi:(1,\infty)\to (1,\infty).$  Now as a direct consequence
of Pisier's characterization of $K$-convex spaces \cite{26} we also have
that an estimate on the $K$-convexity constant of $Y$ in terms of
$T_p(Y)$. Hence we get an estimate of the form
$$ T_2(Y) \le \varphi_p(T_p(Y))C_2(Y^*)$$
for a suitable increasing $\varphi_p:(1,\infty)\to(1,\infty).$  Putting
these estimates together we have
$$ \|T\| \le \psi(T_2(X/E))C_2(Y^*)$$
where $\psi:=\varphi_p\circ\varphi.$
It remains to observe that $C_2(Y^*)\le C_2(Z^*)\le \sqrt2 C_2(X^*)$ and
we are done.\end{proof}

Using this theorem and Pisier's result as in (\ref{e2.9}) we have

\begin{Cor}\label{2.14} $${\Cal E}_X(E;Y)\le
\psi(T_2(X/E))C_2(X^*)C_2(E^*)^{\frac32}C_2(Y)^{\frac32}.$$
\end{Cor}

\section{Linear approximation on convex subsets of $\Bbb
R^n$}\label{linear} \setcounter{equation}0

We begin with the proof of the basic estimate for $w_2(n)$ when $n\ge 2.$
We recall that $w_2(1)=\frac12.$

\begin{Thm}\label{3.1}
We have the estimate:
$$ \frac12\log_2([\frac{n}{2}]+1)\le w_2(n)\le \frac12[\log_2 n]+\frac54.
$$
In particular,
$$ \lim_{n\to\infty}\frac{w_2(n)}{\log_2n}=\frac12.$$
\end{Thm}

{\it Remark.}  See \cite{Ch}, \cite{HU} and \cite{4} for results on the
corresponding problem for convex functions.

In the following discussion $K$ will  denote a closed bounded convex
subset of $\Bbb R^n.$  Note however that our first proposition does not
need convexity:

\begin{Prop}\label{3.2} If $f\in C(K),$ then
$$ E_2(f;K) =\frac12
\max\{\sum_{i=1}^la_if(x_i)-\sum_{j=1}^mb_jf(x_j)\}$$
where the maximum is computed over all pairs of positive integers $l,m$
with $l+m\le n+2$, all subsets $\{x_1,\ldots,x_l\},\ \{y_1,\ldots,y_m\}$
of $K$ and all nonnegative reals $a_1,\ldots,a_l,b_1,\ldots,b_m$ with
$$ \sum_{i=1}^la_i=\sum_{j=1}^mb_j=1 \text{ and }
\sum_{i=1}^la_ix_i=\sum_{j=1}^mb_jy_j.$$
\end{Prop}

\begin{proof}We may choose
 $\varphi$ affine so that
$E_2(f-\varphi;K)=\|f-\varphi\|_K.$  Then clearly $E_2(f-\varphi;K)$
dominates  the expression on the right of the equation.  To prove the
converse, we observe (see, e.g. \cite{27}, p. 36) that there exist
non-empty subsets $\Sigma_+$ and $\Sigma_-$ of $K$ so that
$|\Sigma_+|+|\Sigma_-|\le n+2$ and $(\co \Sigma_+)\cap (\co
\Sigma_-)\neq\emptyset$ and so that for $x\in\Sigma_{\pm}$ we have
$$ f(x)-\varphi(x)=\pm E_2(f;K).$$

Let $\Sigma_+=\{x_1,\ldots,x_l\}$ and $\Sigma_-=\{y_1,\ldots,y_m\}$ then
$l+m\le n+2$  and we can find convex combinations so that
$\sum_{i=1}^la_ix_i=\sum_{j=1}^mb_jy_j.$
Then
\begin{equation}\begin{align*}
E_2(f;K) &=\frac12\{\sum_{i=1}^la_i(f(x_i)-\varphi(x_i))-\sum_{j=1}^mb_j(
f(y_j)-\varphi(y_j))\}\\
&= \frac12\{\sum_{i=1}^la_if(x_i)-\sum_{j=1}^mb_jf(y_j)\}.
\end{align*}\end{equation}
\end{proof}

Let us define $\delta_m:C(K)\to \Bbb R$ (extending the definition of
$\delta_2$) by
\begin{equation}\label{e3.1}
 \delta_m(f):=\sup| f(\sum_{k=1}^ma_kx_k)-\sum_{k=1}^ma_kf(x_k)|
\end{equation}
where the supremum is taken over all $x_1,\ldots,x_m\in K$ and
$a_1,\ldots,a_m\in\Bbb R_+$ such that $\sum_{k=1}^ma_k=1.$
Let
$$ \alpha_m(K)=\sup\{\delta_{m+1}(f): \ f\in C(K), \ \omega_2(f;K)\le
1\}.$$
We then have:

\begin{Cor}\label{3.3}
$$ E_2(f;K)\le \frac12\max\{\alpha_l(K)+\alpha_m(K):\ l,m\ge 0,\
l+m=n\}.$$
\end{Cor}

Observe we have a trivial inequality $\alpha_m(K)\le
\alpha_m(S^m)=:\beta_m$ where $S^m$ is, as usual, the $m$-dimensional
simplex.  Thus, combining with Corollary \ref{3.3} we obtain the
inequality
\begin{equation}
\label{e3.2}
w_2(n)\le \frac12\max\{\beta_l+\beta_m:\ l+m=n\}.
\end{equation}
Note that by Proposition 3.2, $\delta_{m+1}(f)\le 2E_2(f;K)$ and so
$\beta_m\le 2w_2(m).$  In particular, $\beta_1\le 1$ by the results of
Whitney \cite{30}.
To obtain an estimate for all $m$ we need:

\begin{Lem}\label{3.4}
$$ \beta_{2m}\le \beta_m+\frac12, \quad m\in\Bbb N.$$
\end{Lem}

\begin{proof} We set $S^m:=\co\{e_1,\ldots,e_{m+1}\}$ where
$e_1,\ldots,e_{m+1}$ is the canonical basis of $\Bbb R^{m+1}.$  Replacing
$f$ by $f-\varphi$ where $\varphi$ is an affine function satisfying
$\varphi(e_k)=f(e_k)$ for $1\le k\le m+1$ we can obtain an alternate
expression for $\beta_m:$
\begin{equation}\label{e3.3}
\beta_m=\sup\{|f(x)|:\ \omega_2(f)\le 1 \text{ and } f(e_k)=0, \ 1\le
k\le m+1\}.
\end{equation}
Suppose then $f\in C(S^{2m})$ satisfies $\omega_2(f)\le 1$ and $f(e_k)=0$
for $1\le k\le m+1.$  Choose $x\in S^{2m}.$
 Let
$x=\sum_{k=1}^{2m+1}\xi_ke_k.$  Let $(r_k)_{k=1}^{2m+1} $ be a
reordering of $\{1,2,\ldots,2m+1\}$ so that $\xi_{r_k}$ is
increasing.  Then we may choose signs $(\epsilon_k)_{k=1}^{m}$ so that
$$0\le a:=\sum_{k=1}^{m}\epsilon_k(\xi_{r_{2k-1}}-\xi_{r_{2k}})\le
\max_{1\le k\le m}(\xi_{r_{2k}}-\xi_{r_{2k-1}}) $$
Then we can write $x=\frac12(y+z)$ where
$$y= x+\sum_{k=1}^m\epsilon_k(\xi_{r_{2k-1}}e_{r_{2k-1}}
-\xi_{r_{2k}}e_{r_{2k}}) - ae_{r_{2m+1}}$$
and
$$z= x-\sum_{k=1}^m\epsilon_k(\xi_{r_{2k-1}}e_{r_{2k-1}}
-\xi_{r_{2k}}e_{r_{2k}}) + ae_{r_{2m+1}}.$$
Hence
$$|f(x)-\frac12(f(y)+f(z))|\le \frac12.$$
If $y=\sum_{k=1}^{2m+1}\eta_ke_k$ then $\eta_k>0$ at most $m+1$ times
and so $|f(y)-\sum_{k=1}^{2m+1}\eta_kf(e_k)|\le \beta_m.$
With a similar estimate for $z$ we obtain
$$ |f(x)-\sum_{k=1}^{2m+1}\xi_kf(e_k)|\le \beta_m+\frac12.$$
This leads immediately to the claimed estimate.\end{proof}

{\it Proof of the upper estimate in Theorem \ref{3.1}.} Since
$\beta_1=1$, Lemma \ref{3.4} and induction gives us that $\beta_m\le
\frac12\log_2m+1=\frac12k+1$ when $m=2^k.$  Now suppose $2^k\le
n<2^{k+1}$.  Clearly if $l+m=n$ then at most one of them exceeds
$2^{k}.$ Hence
$\beta_l+\beta_m\le \beta_{2^k}+\beta_{2^{k+1}}\le k+\frac52.$  Applying
inequality (\ref{e3.2}) we get the estimate $w_2(n)\le
\frac12[\log_2n]+\frac54.$

For the lower estimate, we require the following general result:

\begin{Lem}\label{3.6} Suppose $K_1,K_2$ are closed bounded convex
subsets of $\Bbb R^{n_1},\Bbb R^{n_2},$ respectively. Suppose $f_i\in
C(K_i)$ for
$i=1,2$ are convex
and satisfy $\omega_2(f_i;K_i)\le 1.$  Then if $g:K_1\times K_2\to \Bbb
R$ is defined by $g(x,y)=f_1(x)-f_2(y)$ we have:\newline
(a) $E_2(g;K_1\times K_2) = E_2(f_1;K_1)+E_2(f_2;K_2)$\newline
(b) $\omega_2(g;K_1\times K_2)\le 1.$
\end{Lem}

\begin{proof} Suppose $h=(h_1,h_2)\in \Bbb R^{n_1}\times\Bbb R^{n_2}.$
Then
$$ \Delta^2_hg(x,y)=\Delta_{h_1}^2f_1(x)-\Delta_{h_2}^2f_2(y).$$
Since $f_1$, $f_2$ are convex we obtain $\omega_2(g;K_1\times K_2)\le 1.$
This proves (b).  To prove (a) it suffices to apply Proposition \ref{3.2}
(cf. Theorem 6.2.5 in \cite{27}).\end{proof}

{\it Proof of the lower estimate in Theorem \ref{3.1}.} Let
$S^n=\co\{e_1,\ldots,e_{n+1}\}$ as before.  Define the function
\begin{equation}\label{e3.7}
f_n(x):=\frac12\Lambda(x)=\frac12\sum_{k=1}^{n+1}x_k\log_2x_k.\end{equation}
Since the function $\psi(t):=t\log_2t,\ 0\le t\le 1,$ satisfies $0\le
\Delta_h^2f(t)\le \Delta_{|h|}^2f(0)=2|h|$, the function $f_n$ is convex
and
$$ 0\le \Delta_h^2f_n(x)=\frac12\sum_{k=1}^{n+1}\Delta_{h_k}^2
\psi(x_k)\le \sum_{k=1}^{n+1}|h_k|.$$
Now $h=\frac12((x+2h)-x)$ so that $\sum_{k=1}^{n+1}|h_k|\le 1.$  Thus
$\omega_2(f_n)\le 1.$  Let $u:=\frac{1}{n+1}\sum_{k=1}^{n+1}e_k.$  Then
by Proposition \ref{3.2},
$$ E_2(f_n) \ge \frac{1}{n+1}\sum_{k=1}^{n+1}f_n(e_k)-
f_n(u)=\frac14\log_2(n+1).$$  We remark that this function was
essentially first considered in this context (in an equivalent
formulation) by Ribe \cite{Rib}.

We can now apply Lemma \ref{3.6}.  If $n=2m$, putting $K_1=K_2=S^m$ and
using $f_n$ for both $f_1$ and $f_2$ of the Lemma,
we obtain the existence of $g$ on $S^n\times S^n$ with $\omega_2(g)\le
1$ and $E_2(g)\ge \frac12\log_2(m+1).$  Hence $\omega_2(n)\ge
\frac12\log_2(\frac{n}{2}+1).$  If $n=2m+1$ we put $K_1=S^m,$ and
$K_2=S^{m+1}$ and use $f_m, f_{m+1}$ to deduce that
$$ \omega_2(n)\ge \frac14(\log_2(m+1)+\log_2(m+2))\ge
\frac12\log_2([\frac{n}{2}]+1).$$
The proof of Theorem \ref{3.1} is now complete. \qed
\smallskip

{\bf Remarks.} (a) For small values of $n$ we can use (\ref{e3.2})
directly to obtain better upper bounds for $w_2(n).$  Thus $\beta_2\le
\frac32,\ \beta_3,\beta_4\le 2$ and hence $w_2(2)\le 1$, $w_2(3)\le
\frac54$ and $w_2(4)\le \frac32.$

 On the other hand, if we use the
piecewise linear function
$f_{\epsilon}(t)=\max((1-\epsilon)(1-\epsilon^{-1}t),0)$ on
$[0,1]$ then $f_{\epsilon}$ is convex and satisfies
$\omega_2(f_{\epsilon})=1$ and $E_2(f_{\epsilon})=\frac12(1-\epsilon).$
Then using Lemma \ref{3.6} and the functions
$g_{\epsilon}(x,y)=f_{\epsilon}(x)-f_{\epsilon}(y)$ we obtain that
$w_2([0,1]^2)\ge 1-\epsilon.$  Combined with the upper estimates above we
obtain:
\begin{equation}\label{e3.8}
w_2(2)=w_2([0,1]^2)=1.
\end{equation}
Notice that $w_2([0,1]^2)=w_2(\ell_{\infty}^2)=w_2(\ell_1^2).$

(b) The corresponding examples considered in \cite{7} show that
$\beta_2=\frac53$ and $\beta_3=2.$

We now show that for the case of the simplex the lower bound
$\frac14\log_2(n+1)$ is asymptotically sharp.  More precisely:

\begin{Thm}\label{3.7} $$
\lim_{n\to\infty}\frac{w_2(S^n)}{\log_2n}=\frac14.$$
\end{Thm}

We remark first that the functions $f_n$ constructed in (\ref{e3.7}) show
that $w_2(S^n)\ge \frac14\log_2(n+1)$ so that
\begin{equation}\label{e3.9}
\liminf_{n\to\infty}\frac{w_2(S^n)}{\log_2n} \ge \frac14.
\end{equation}
The proof of Theorem \ref{3.7} will follow from the following Theorem
\ref{2.9}.

\begin{Thm}\label{2.9} For any $0<\epsilon<\frac12$ there is a
constant $C=C(\epsilon)$ such that
 whenever $n\in\Bbb N,$ and $f\in C(S^n)$
satisfies
$\omega_2(f)\le 1-\epsilon$ there is a Banach lattice $X$ so that
$$ |f(u)-\frac12(\Phi_X(u)-\Phi_{X^*}(u))|\le C $$
for all $u\in S^n.$
\end{Thm}

 Before proving Theorem \ref{2.9} let us complete the proof of Theorem
\ref{3.7} assuming Theorem \ref{2.9}:

{\it Proof of Theorem of \ref{3.7}:} Fix $\epsilon>0.$  If $f\in C(S^n)$
satisfies $\omega_2(f)\le 1-\epsilon,$ we determine $X$ so that Theorem
\ref{2.9} holds.  Let $\|\cdot\|_E$ be the Hilbertian norm determined by
the John ellipsoid for $B_X$ \cite{8}. Then in the terminology of
Proposition
\ref{2.11} we must have $E=\ell_2(a)$ for a suitable positive sequence
$a=(a_1,\ldots,a_{n+1}).$  Then by Proposition \ref{2.11} we have that
$\Phi_E-\Phi_{E^*}$ is linear: in fact $\Phi_E(u)-\Phi_{E^*}(u)=-2\langle
u,\log a\rangle.$

From the properties of the John ellipsoid we have $B_E\subset B_X\subset
(n+1)^{1/2}B_E$ so that $\Phi_E(u)\le \Phi_X(u)\le
\Phi_E(u)+\frac12\log_2(n+1).$  From Theorem \ref{2.10} (a) we
get
$$\Phi_{E^*}(u)-\frac12 \log_2(n+1) \le \Phi_{X^*}(u) \le \Phi_{E^*}(u)$$
and
so
$$ |\frac12(\Phi_X(u)-\Phi_{X^*}(u))+\langle u,\log a\rangle|\le
\frac14\log_2(n+1).$$
It follows that
$$ E_2(f) \le C(\epsilon)+\frac14\log_2(n+1).$$
This implies that
$$w_2(S^n)\le (1-\epsilon)^{-1}(C(\epsilon)+\frac14\log_2(n+1),$$
which in turn gives the required upper estimate
$$ \limsup_{n\to\infty}  \frac{w_2(S^n)}{\log_2n}\le \frac14.$$
This completes the proof of Theorem \ref{3.7}.\qed

We now turn to the proof of Theorem \ref{2.9}.

\begin{proof} Let $f:S^n\to \Bbb R$ be a bounded function satisfying the
condition $\omega_2(f)\le 1-\epsilon$ where $0<\epsilon<\frac12$  is
fixed.
By Whitney's theorem applied to each line segment we have $\delta_2(f)\le
\omega_2(f)\le 1-\epsilon.$  Let $\alpha:=1-\frac12\epsilon$ and apply
Theorem \ref {2.8} to the function $\alpha^{-1}f$.  Thus there is an
$(n+1)-$dimensional Banach lattice $Y$ with
\begin{equation}\label{e3.12}
\|f-\alpha(\Phi_Y-\Phi_{Y^*})\|_{S^n}\le C(\epsilon).
\end{equation}
To complete the proof we will find a lattice $X$ for which
\begin{equation}\label{e3.13}
\|\alpha(\Phi_Y-\Phi_{Y^*})-\frac12(\Phi_X-\Phi_{X^*})\|_{S^n}\le
C'(\epsilon).
\end{equation}

In order to do this we will show the existence of a Banach lattice $X$
such that if we put $\theta:=1-(2\alpha)^{-1}$ then the spaces $Y$ and
$X^{1-\theta}(\ell_2^{n+1})^{\theta}$ have equivalent norms with the
constant of equivalence depending only on $\epsilon.$  Assuming this
fact, let us show how the proof is completed.  In this case by Theorem
\ref{2.10} and Proposition \ref{2.11} we have
$$ \|\Phi_Y-((1-\theta)\Phi_X+\frac{\theta}{2}\Lambda)\|_{S^n}\le
C_1(\epsilon).$$  Using the duality result Theorem \ref{2.10} (a)
this implies that
$$ \|(\Phi_Y-\Phi_{Y^*})-(1-\theta)(\Phi_X-\Phi_{X^*})\|_{S^n}\le
2C_1(\epsilon).$$
Since $1-\theta=(2\alpha)^{-1}$ this establishes (\ref{e3.13}) and
combined with (\ref{e3.12}) the theorem is proved.

Thus it remains to construct $X.$
We will need the following lemma:

\begin{Lem}\label{L3.8} Suppose $p$ is defined by
$$p:=(1+\frac{1}{2\alpha})^{-1}+(1+\frac{1-\epsilon}{2\alpha})^{-1}.$$
 Then there is a constant $C$ depending only on $\epsilon$ so that for
every disjoint family of vectors $\{y_i\}_{i=1}^m\subset \Bbb R^{n+1},$
we have
\begin{equation}\label{e3.17}
\|\sum_{i=1}^my_i\|_Y \le C(\sum_{i=1}^m\|y_i\|_Y^p)^{1/p}
\end{equation}
and
\begin{equation}\label{e3.17b}
\|\sum_{i=1}^my_i\|_{Y^*} \le C(\sum_{i=1}^m\|y_i\|_{Y^*}^p)^{1/p}
\end{equation}
\end{Lem}

Before proving the lemma, let us show how to complete the construction
of
$X$ assuming this lemma.  We set
\begin{equation}\label{e3.15}
\frac1r:=\frac12+\frac1{4\alpha}=1-\frac{\theta}{2}
\end{equation}
Then $p>r$.  By the lemma both $Y$ and $Y^*$ satisfy upper $p$-estimates
with constants depending only on $\epsilon.$  According to a well-known
theorem of Maurey and Pisier (see, e.g. \cite{17}) this implies that $Y$
and $Y^*$ are both $r$-convex with constants depending only $\epsilon.$
This means that for any $y_1,\ldots,y_m\in Y$ we have
$$ \|(\sum_{i=1}^m|y_i|^r)^{\frac{1}{r}}\|_Y \le C
(\sum_{i=1}^m\|y_i\|_Y^r)^{1/r}$$
where $C$ depends only on $\epsilon$, and a similar inequality holds in
$Y^*.$  Now by Propositions 1.d.4 and 1.d.8 of \cite{17} there is
a lattice $Y_0$ so that $Y_0,Y_0^*$ are $r$-convex with constant one and
the $Y_0$-norm is $C$-equivalent to the $Y$-norm with $C$ depending only
on $\epsilon.$  Finally we use the Pisier extrapolation theorem \cite{23}
to deduce that there is a Banach lattice $X$ so that
$Y_0=X^{1-\theta}(\ell_2^{n+1})^{\theta}.$\end{proof}

We now turn to the proof of the Lemma.

\begin{proof} Let $g:=\Phi_Y-\Phi_{Y^*}.$  Using (\ref{e3.12}) we first
estimate $\delta_m(g)\le \alpha^{-1}(C+\delta_m(f))$ where $C$ depends
only on $\epsilon.$ Since $\omega_2(f)\le 1-\epsilon$, Lemma \ref{3.4}
gives that $\delta_m(f)\le (1-\epsilon)(\frac12\log_2m+1).$  Hence
\begin{equation}
\label{e3.20}
\delta_m(g)\le C_1+\frac{1-\epsilon}{2-\epsilon}\log_2m
\end{equation}
where $C_1=C_1(\epsilon).$

Now suppose $u_1,u_2,\ldots,u_m\in S^n$ have disjoint supports and that
$u=\sum_{i=1}^ma_iu_i\in S^n$ is a convex combination.  Then
$$ g(u) \ge
\sum_{i=1}^ma_ig(u_i)-C_1-\frac{1-\epsilon}{2-\epsilon}\log_2m.$$
By duality (Proposition \ref{2.10}) $\Phi_Y=\frac12(g+\Lambda)$ and
direct calculation gives  us that:
\begin{equation}\begin{align*}
 \Lambda(u) &= \sum_{i=1}^ma_i\Lambda(u_i)+\sum_{i=1}^ma_i\log_2a_i
\quad\quad\quad\quad\\
 &\ge
\sum_{i=1}^ma_i\Lambda(u_i)-\log_2m.
\end{align*}\end{equation}
Combining these estimates we have
\begin{equation}\label{e3.21}
\Phi_Y(u)\ge \sum_{i=1}^ma_i\Phi_Y(u_i)-\frac12C_1-\frac1{p_1}\log_2m
\end{equation}
where
\begin{equation}\label{e3.22}
\frac{1}{p_1}:=\frac12+\frac{1-\epsilon}{2(2-\epsilon)}\end{equation}
Note that we have a precisely similar estimate to (\ref{e3.21}) for
$\Phi_{Y^*}$ in place of $\Phi_Y$  using instead the equation
$\Phi_{Y^*}=\frac12(\Lambda-g).$

Now suppose $y_1,\ldots,y_m\in Y$ have disjoint supports.  For any $u\in
S^n$ we can write $u=\sum_{i=1}^ma_iu_i$ as a convex combination where
 supp $u_i\supset$ supp $y_i$ and
the $(u_i)$ have disjoint supports.  Let
us write $$ \langle v,\log_2|x|\rangle:=\sum_{i=1}^{n+1}v_i\log_2|x_i|$$
for $v\in S^n$ and $x\in R^{n+1}$ (with $-\infty$ as a possible value!).
 Then by (\ref{e3.21})
\begin{equation}\begin{align*}
 \Phi_Y(u) &\ge \sum_{i=1}^ma_i\langle u_i,\log_2|y_i|\rangle
+\log_2(2^{-C_2}m^{-\frac{1}{p_1}})\\
&=\langle
u,\log_2(2^{-C_2}m^{-\frac{1}{p_1}}|y_1+\cdots+y_m|)\rangle
\end{align*}\end{equation}
where $C_2:=\frac12C_1.$
Now it is a consequence of Theorem 4.4 of \cite{11} (which is much
simpler in our finite-dimensional setting) that this implies
\begin{equation}\label{e3.23}
\|y_1+\cdots+y_m\|_Y \le 2^{C_2}m^{\frac{1}{p_1}}
\end{equation}
Again the same inequality holds in $Y^*.$

Now suppose $\{y_1,\ldots,y_m\}$ are any disjoint vectors with
$\sum_{i=1}^m\|y_i\|^p_Y= 1.$
Let $A_k:=\{i:2^{-k}<\|y_i\|_Y\le 2^{1-k}\}.$  If $|A_k|$ denotes the
cardinality of $A_k$ then $|A_k|\le 2^{kp}$ and by (\ref{e3.21}) we have
$$\|\sum_{i=1}^my_i\|_Y \le \sum_{k=1}^{\infty}\|\sum_{i\in
A_k}y_i\|_Y
\le
  \sum_{k=1}^{\infty}2^{C_2+1}|A_k|^{\frac{1}{p_1}}2^{-k}
 \le 2^{C_2+1}\sum_{k=1}^{\infty}2^{k(\frac{p}{p_1}-1)}.
$$
Since
$$\frac{p}{p_1}=\frac12(\frac{3-2\epsilon}{3-\epsilon}+1)<1$$ this
implies an estimate
$$ \|\sum_{i=1}^my_i\|\le C(\epsilon)<\infty$$
and, combined with the similar estimate for $Y^*$, this establishes the
lemma.
\end{proof}

We now turn our attention to the case when $K=B_X$ is the unit ball of a
finite-dimensional Banach space.  Our main result concerns the case when
$X=\ell_p^n$ for $1\le p\le \infty.$

\begin{Thm}\label{3.8} (a) There exist constants $c_1,c_2>0$ so that
$$ c_1\log_2(n+1)\le w_2(\ell_1^n)\le c_2\log_2(n+1).$$
In addition
\begin{equation}\label{e3.24}
\limsup_{n\to\infty}\frac{w_2(\ell_1^n)}{\log_2n}\le \frac14
\end{equation}
(b) If $1<p\le 2$ then $\gamma(p):=\sup_{n\in\Bbb N}w_2(\ell_p^n)<\infty$
and further there exist constants $d_1,d_2>0$ so that
$$ \frac{d_1}{p-1}\le \gamma(p)\le \frac{d_2}{p-1}|\log(p-1)|.$$
(c) If $2\le p\le\infty$ then $\gamma(p):=\sup_{n\in\Bbb
N}w_2(\ell_p^n)<\infty$, and further
$$ \gamma:=\sup_{2\le p\le\infty}\gamma(p)\le 1602<\infty.$$
\end{Thm}

{\it Proof of (a).} The upper estimate is an immediate consequence of
Theorem \ref{3.1}.  To prove the lower estimate, let $\tilde
f_n:B_{\ell_1^n}\to\Bbb R$ be defined by
$$ \tilde f_n(x):=\frac12
\sum_{i=1}^nx_i\log_2|x_i|=\frac12\sum_{i=1}\psi(t)$$
where $\psi(t):= t\log_2|t|$ for $-1\le t\le 1.$  Since
$|\Delta_h^2\psi(t)|\le 2\log_2(1+\sqrt2)|h|$, we obtain
$$|\Delta_h^2\tilde f_n(x)|\le \log_2(1+\sqrt 2)\sum_{i=1}|h_i|$$
so that $\omega_2(\tilde f_n)\le \log_2(1+\sqrt2 ).$  Since $\tilde
f_n|_{S^{n-1}}=f_{n-1}$ as defined in (\ref{e3.7}) we have
$$ E_2(\tilde f_n;B_{\ell_1})\ge E_2(f_{n-1};S^{n-1})\ge
\frac14\log_2n.$$
This implies that
$$ w_2(\ell_1^n)\ge \frac{1}{4\log_2(1+\sqrt 2)}\log_2n.$$

It remains to prove (\ref{e3.24}).  For this we need:

\begin{Lem}\label{3.9} $w_2(\ell_1^n)\le w_2(S^{n-1})+\frac32.$\end{Lem}

\begin{proof}
Suppose first $f$ is a bounded continuous function on $B_{\ell_1^n}$
with $\omega_2(f)\le 1.$   Then there is an affine function $g$ defined
on $S^{n-1}$ with $|\frac12(f(x)-f(-x))-g(x)| \le w_2(S^{n-1}).$
We can extend $g$ to a linear functional on $\ell_1^n.$
We also have $|\frac12(f(x)+f(-x))-f(0)|\le \frac12$ for $x\in S^n.$
Hence $|f(x)-g(x)-f(0)|\le w_2(S^{n-1})+\frac12$ if $x\in \pm S^{n-1}.$
Now if $x\in B_{\ell_1^n}$ we can find $u,v\in S^{n-1}$ and $0\le t\le 1$
so that $x=tu-(1-t)v.$  Hence $|f(x)-tf(u)-(1-t)f(-v)|\le 1$ by the
one-dimensional Whitney result which is essentially the fact that
$\delta_2(f)\le\beta_1=1.$ Since
$g$ is linear,
$|f(x)-g(x)-f(0)|\le
w_2(S^{n-1})+\frac32.$  It follows that $w_2(\ell_1^n)\le
w_2(S^{n-1})+\frac32.$ \end{proof}

Now the inequality (\ref{e3.24}) follows from Theorem \ref{3.7} and the
proof of (a) is complete. \qed

We postpone the proof of (b) until after (c):

{\it Proof of (c):}We will need the following Lemma (see (\ref{e2.1})
for the definition of $v_2(X)$):

\begin{Lem}\label{3.10}
(a) Let $E$ be a subspace of a finite-dimensional Banach space $X.$  Then
$v_2(X/E)\le 2v_2(X).$ \newline
(b) Suppose $X,Y$ are two $n$-dimensional  Banach spaces.  Then
$v_2(Y)\le d(X,Y)v_2(X)$ where $d(X,Y)$ is the Banach-Mazur distance
between $X$ and $Y.$
\end{Lem}

\begin{proof} Let $Q:X\to X/E$ be the quotient map.  If $f:X/E\to\Bbb R$
is a continuous homogeneous function then there is a linear functional
$x^*$ on $X$ so that
$$ |f(Qx)-x^*(x)|\le v_2(X)\omega_2(f;B_{X/E})\|x\|$$
for $x\in X.$  For $x\in E$ we have
$$ |x^*(x)|\le v_2(X)\omega_2(f)\|x\|$$ and so by the Hahn-Banach theorem
we can find a linear functional $u^*$ with $u^*(e)=x^*(e)$ for $e\in E$
and $\|u^*\|\le v_2(X)\omega_2(f).$  Then there exists $z^*\in (X/E)^*$
with $x^*-u^*=z^*\circ Q$ and we have:
$$|f(Qx)-z^*(Qx)|\le  |f(Qx)-x^*(x)|+|u^*(x)|\le
2v_2(X)\omega_2(f)\|x\|.$$
Part (a) now follows.

For part (b) suppose $T:X\to Y$ satisfies $\|T\|=1$ and
$\|T^{-1}\|=d(X,Y).$  Then if $f:Y\to\Bbb R$ is a continuous homogeneous
function then $\omega_2(f\circ T;B_X)\le \omega_2(f;B_Y).$  Now there
exists $x^*\in X^*$ so that $|f(Tx)-x^*(x)|\le
v_2(X)\omega_2(f;B_Y)\|x\|.$
Let $y^*=x^*\circ T^{-1}.$  Then $|f(y)-y^*(y)| \le
v_2(X)\omega_2(f;B_Y)d(X,Y)\|y\|$ and the lemma follows.
\qed

Now suppose $2\le p\le\infty.$ Then for any $n\in\Bbb N$ and
$\epsilon>0$ there exists $N$ so that $\ell_p^n$ is
$(1+\epsilon)-$isomorphic to a quotient of $\ell_{\infty}^N$.  Hence
$v_2(\ell_p^n)\le 2(1+\epsilon)v_2(\ell_{\infty}^N).$

However, the estimate $v_2(\ell_{\infty}^N)\le 200$ is proved in
\cite{14}
(a factor 2 was omitted from the argument as pointed out in \cite{18}).
Hence $v_2(\ell_p^n)\le 400$ for all $n.$  Now by Proposition \ref{2.1}
we have
$$ w_2(\ell_p^n)\le 4v_2(\ell_p^n)+\frac32\le 1602.$$  (Note that for
$p=\infty$ we can eliminate a factor of 2 and get an estimate of 802.)
\end{proof}

We now proceed to the proof of (b).  Let us comment first that there is a
striking difference between the cases $p<2$ and $p>2$ and this reflects
the differing behavior of these spaces with respect to (Rademacher) type
(see Section 2 for the definitions.)

We start by establishing the lower bound.  For this we note that
$d(\ell_1^n,\ell_p^n)\le n^{1/q}$ where $\frac1p+\frac1q=1.$  Hence
by
Lemma \ref{3.10} and part (a) we have $v_2(\ell_p^n)\ge
c_1n^{-1/q}\log(n+1).$
  If we choose $n=[e^q]$ we obtain an estimate $\gamma(p)\ge d_1q\ge
d_1(p-1)^{-1}$ where $d_1>0.$

We will derive the upper bound from a general result about the
relationship between the Whitney constants and the Rademacher type $p$
constant.

\begin{Thm}\label{3.11} There is an absolute constant $C$ so that if
$X$ be a finite-dimensional Banach space and  $1<p\le 2,$ then
\begin{equation}\label{3.27}
w_2(X) \le \frac{C}{p-1}(1+|\log(p-1)|+\log T_p(X))
\end{equation}
\end{Thm}

\begin{proof} For this theorem we need the following elementary lemma:

\begin{Lem}\label{type}
Suppose $Y$ is a Banach space of type $p$ where $1<p\le 2$ with type $p$
constant $T_p(Y).$  Suppose $y_1,\ldots,y_n\in B_Y$ and that $k\in\Bbb
N.$
Then there is a subset $\sigma$ of $\{1,2,\ldots,n\}$ with $|\sigma|\le
2^{-k}n$ and so that
$$ \|\frac1n\sum_{i=1}^ny_i-\frac{2^k}{n}\sum_{i\in\sigma}y_i\|\le
T_p(Y)n^{-1/q}\frac{2^{k/q}-1}{2^{1/q}-1},$$ where, as usual,
$\frac1p+\frac1q=1.$
\end{Lem}

\begin{proof} We prove this by induction on $k,$ with $k=0$ as the
trivial starting point.  Suppose $\sigma_k$ is the subset satisfying the
conclusions of the lemma for $k.$  Then by the definition of the type
$p$ constant there is a choice of signs $\epsilon_i=\pm1$ with
$$ \|\sum_{i\in\sigma_k}\epsilon_i y_i\| \le T_p(Y)|\sigma_k|^{1/p}\le
T_p(Y)2^{-k/p}n^{1/p}.$$
Without loss of generality we can assume
$\sum_{i\in\sigma_k}\epsilon_i\le
0.$  Let $\sigma_{k+1}:=\{i\in\sigma:\epsilon_i=1\}.$ Then
\begin{equation*}
\|\frac{2^k}{n}\sum_{i\in\sigma_k}y_i-\frac{2^{k+1}}{n}\sum_{i\in\sigma_
{k+1}}y_i\| = \frac{2^k}{n}\|\sum_{i\in\sigma_k}\epsilon_iy_i\|
\le T_p(Y) 2^{k/q}n^{-1/q}.
\end{equation*}
The induction step now follows easily\end{proof}

Returning the proof of Theorem \ref{3.11}, we will estimate
$v_2:=v_2(X).$  Suppose that $f$ is any continuous homogeneous function
on $X$ with $\omega_2(f;B_X)\le 1.$
We may pick $x^*\in X^*$ so that
if $g:=f-x^*,$ then
\begin{equation}\label{ee3.18}
E_2(f;B_X)=E_2(g:B_X)=\sup\{|g(x)|:\|x\|\le 1\}\le
v_2.\end{equation}

By Proposition \ref{3.2},
$$ E_2(g;B_X)\le \sup\{\delta_m(g;B_X):\ m\in\Bbb N\}$$ where
$\delta_m(f;B_X)$ is defined in (\ref{3.1}).   Since $g$ is continuous
the right-hand side is equal to $\sup_{n\in \Bbb N}b_n$  where
$$ b_n:=\sup\{|g(\frac1m\sum_{i=1}^mx_i)-\frac1m\sum_{i=1}^mg(x_i)| : \
x_1,\ldots,x_m\in B_X,\ m\le n\}.$$

We will show that
\begin{equation}\label{ee3.19}
b_n\le 3+40q +2q\log T_p +2q\log v_2\end{equation}
where $T_p:=T_p(X).$

To establish (\ref{ee3.19}) choose an integer $N:=[(T_pv_2)^q].$  By
Theorem \ref{3.1} $b_n\le 2w_2(n)\le 3+\log_2n$ and this shows that
 $$b_N \le 3+q\log_2T_p +q\log_2v_2\le 3+2q\log T_p+2q\log v_2.$$
In particular, (\ref{ee3.19}) holds for all $n\le N.$

Suppose  now $n>N$ and choose $k\in\Bbb N$ so that $2^{k-1}N<n\le
2^kN.$
We consider the space $Y:=X\oplus_{\infty}\Bbb R$; then $T_p(Y)\le
2T_p(X)=2T_p.$  If $x_1,\ldots,x_n\in B_X$ we define elements of $B_Y$ by
$y_i:=(x_i,v_2^{-1}g(x_i)).$  By Lemma \ref{type} there is a
subset
$\sigma$ of $\{1,2,\ldots,n\}$ with $|\sigma|\le 2^{-k}n$ so that
\begin{equation}\label{ee3.20}
 \|\frac1n\sum_{i=1}^ny_i-\frac{2^k}{n}\sum_{i\in\sigma}y_i\|_Y \le
q (\log 2)^{-1}T_p(Y)2^{k/q}n^{-1/q}
\le 8qT_pN^{-1/q}.\end{equation}
In particular, we have if $u:=\frac1n\sum_{i=1}^nx_i$ and
$w:=\frac{2^k}{n}\sum_{i\in\sigma}x_i$
\begin{equation}\label{eq1} \|u-w\|\le
 8qT_pN^{-1/q}.\end{equation}
Since
$u,w\in B_X$ and $g$ is homogeneous, we have
$|g(u-w)-g(u)+g(w)|\le \omega_2(f;B_X)\le
1.$  Hence and by (\ref{ee3.18}),
\begin{equation}\label{eq2}|g(u)-g(w)| \le 1+
|g(u-w)| \le 1+v_2\|u-w\|\le 20q
\end{equation} by the choice of $N.$
 We also have  from (\ref{ee3.20})
\begin{equation}\label{eq3}
|\frac1n\sum_{i=1}^ng(x_i)-\frac{2^k}{n}\sum_{i\in\sigma} g(x_i)|\le
8qv_2T_pN^{-1/q} \le 20q.\end{equation}
Finally we note that, since $|\sigma|\le 2^{-k}n\le N,$
\begin{equation}\label{eq4}
|g(w)-\frac{2^k}{n}\sum_{i\in\sigma}g(x_i)|\le b_N\le 3+2q\log
T_p+2q\log v_2\end{equation}

Combining (\ref{eq2}), (\ref{eq3}) and (\ref{eq4}) gives us
\begin{equation*}
|g(u)-\frac1n\sum_{i=1}^ng(x_i)|\le 3+40q  +2q\log T_p+2q\log
v_2\end{equation*} and so (\ref{ee3.19}) holds.

Now (\ref{ee3.19}) gives an estimate independent of $n$ and so
implies that
$$ E_2(f;B_X)=E_2(g:B_X)\le
3+ 40q  +2q\log T_p+2q\log
v_2.$$  Since this estimate holds for all such $f,$ we obtain
$$ v_2 \le 3+40q + 2q\log T_p +2q\log v_2.$$
Since $q\log v_2\le \frac14v_2+q\log q+q\log 4$ this gives the required
upper estimate in (b).
\end{proof}

The proof of Theorem \ref{3.11} is now also complete.\qed

\begin{Cor}\label{3.15}
There is a universal constant $C$ so that if $X$ is an $n$-dimensional
Banach space  and $2<q<\infty,$
$$ w_2(X)\le Cq(\log q+\log C_q(X^*)+\log(1+\log n)).$$
\end{Cor}

\begin{proof} If $\frac1p+\frac1q=1$ then $T_p(X)\le C(\log n+1)C_q(X^*)$
(see \cite{28}).  It remains to apply the inequality
(\ref{3.27}).\end{proof}

Note that for the case of $\ell_{\infty}^n$ this is weaker than the
conclusion of Theorem \ref{3.8} (c).  We conjecture that there is an
estimate of the form $w_2(X)\le \phi(q,C_q(X^*))$ for a suitable function
$\phi.$  It is possible that the estimate $w_2(X)\le Cq(1+\log C_q(X^*))$
holds, which  would imply  $w_2(X)\le C(p-1)^{-1}(1+\log T_p(X))$ and
$w_2(\ell_p^n)\le C(p-1)^{-1}$ giving a sharp estimate for
$w_2(\ell_p^n).$

 \section{Quadratic approximation on symmetric convex bodies}
\setcounter{equation}0

We now consider the problem of estimating $w_3(X)$ when $X$ is a
finite-dimensional Banach space.  Our first result gives a quite sharp
estimate of $w_3^{(s)}(n):=\sup_{\dim X=n}w_3(X).$

\begin{Thm}\label{4.1} There are absolute constants $0<c,C<\infty$ so
that for every $n\ge 1$
$$ c\sqrt n \le w_3^{(s)}(n)\le C\sqrt n \log (n+1).$$
\end{Thm}

\begin{proof} The upper estimate is a special case of Theorem
\ref{5.2}, (or Corollary \ref{5.6})
which we therefore postpone to the next section. For the lower estimate,
we  use the fact that the space $\ell_1^n$ contains a subspace $V$ so that
every linear projection $P:\ell_1^{n}\to V$
satisfies
\begin{equation}\label{e4.1}
\|P\|\ge c\sqrt n
\end{equation}
where $c>0$ is an absolute constant. This follows from a well-known
result of Kashin \cite{15} that we may pick $V$ with $\dim
V=[\frac{n}{2}]$ and $d(V,\ell_2^{\dim V})\le C$ where $C$ is independent
of
$n.$ For convenience let
$Y$ be the space
$\Bbb R^{n}$ with the norm, 2-equivalent to the $\ell_1-$norm,
$$ \|x\|_Y:=\|x\|_{\ell_1^{2n}}+ \|x\|_{\ell_2^{2n}}.$$
Then (\ref{e4.1}) holds for every linear projection $P:Y\to V$, with
perhaps a different constant.  Since $Y$ is strictly convex, for every
$x\in \Bbb R^n$ there is a unique $\Omega(x)\in V$ so that
$$ \|x-\Omega(x)\|_Y =d_Y(x, V):=\inf_{v\in V}\|x-v\|_Y.$$
The map $\Omega$ is called the {\it metric projection} of $Y$ onto $V$
and
the following properties are well-known (see, e.g. \cite{27}, Sec 5.1):

\begin{Lem}\label{4.2} (a) $\Omega$ is homogeneous and
continuous;\newline
(b) $\Omega$ is a (nonlinear) projection, $\|\Omega(x)\|_Y\le 2\|x\|_Y$
for
$x\in Y$ and $\Omega(x+v)=\Omega(x)+v$ if $x\in Y,\ v\in V.$\newline
(c) For $x,y\in Y,$
$$ \|\Omega(x+y)-\Omega(x)-\Omega(y)\|\le 2(d_Y(x,V)+d_Y(y,V)).$$
\end{Lem}

Now let $\langle,\rangle$ be the standard inner-product on $\Bbb R^n.$
Let $\pi$ be the orthogonal projection onto $V$ and let $\pi^{\perp}$ be
the complementary projection onto $V^{\perp}.$  Let
$\|x\|_{Y^*}:=\sup\{\langle x,y\rangle:\ \|y\|_Y\le 1\}$ be the dual
norm on $\Bbb R^n.$

We now define a norm $\|\ \|_X$ on $\Bbb R^n$ by the formula:
\begin{equation}\label{enorm}
\|x\|_X := d_{Y^*}(\pi x,V^{\perp})+d_{Y}(\pi^{\perp}x,V)
\end{equation}
where
$$ d_{Y^*}(x,V^{\perp})=\inf\{\|x-v^{\perp}\|_{Y^*}:\ v^{\perp}\in
V^{\perp}\}.$$

Finally let us define the continuous homogeneous function
\begin{equation}\label{edef}
F(x) := \langle \pi x ,\Omega(\pi^{\perp} x ) \rangle.
\end{equation}

Now suppose $x,x+3h\in B_X.$  Let $x=x_1+x_2$ and $h=h_1+h_2$ where
$x_1,h_1\in V$ and $x_2,h_2\in V^{\perp}.$  Then
\begin{equation}\label{4.a}
 \Delta_h^3F(x)=\langle x_1,\Delta_{h_2}^3\Omega(x_2)\rangle + 3\langle
h_1,\Delta_{h_2}^2\Omega(x_2+h_2)\rangle.\end{equation}
Now we have
\begin{equation}
\begin{align*}
 \|\Delta_{h_2}^3\Omega(x_2)\|_Y &\le \|\Delta_{h_2}^2\Omega(x_2)\|_Y +
\|\Delta_{h_2}^2\Omega(x_2+h_2)\|_Y\\
&\le 2(d_Y(x_2,V)+d_Y(x_2+2h_2,V)+d_Y(x_2+h_2,V) +d_Y(x_2+3h_2,V))\\
&\le 8 \end{align*}
\end{equation}
by Lemma \ref{4.2}. Similarly
$$ \|\Delta_{h_2}^2\Omega(x_2+h_2)\|\le 4.$$
Hence by (\ref {4.a}) have
\begin{equation}\label {4.b}
|\Delta_h^3F(x)| \le 8d_{Y^*}(x_1,V^{\perp})+ 12d_{Y^*}(h_1,V^{\perp})\le
16\end{equation}
since $d_{Y^*}(x_1,V^{\perp})\le 1$ and $d_{Y^*}(h_1,V^{\perp})\le
\frac{2}{3}.$
Thus (\ref {4.b}) implies
\begin{equation}\label{4.c}
\omega_3(F;B_X)\le 16.
\end{equation}
Let $v_3:=v_3(X).$  Then there is a quadratic form $Q(x)$ such that
$$ |F(x)-Q(x)| \le 16v_3\|x\|_X^2$$
for $x\in \Bbb R^n.$  We can write $Q(x)=\langle x,Ax\rangle$ where $A$
is a symmetric $n\times n$ matrix or equivalently a symmetric linear
operator on $\Bbb R^n.$

Note for every $x\in \Bbb R^n$ we have $F(\pi x )=F(\pi^{\perp} x )=0.$
Hence
$$ |\langle \pi x , A\pi x\rangle|\le 16v_3\|\pi x \|_X^2 \le
16v_3\|x\|_X^2$$
and
$$ |\langle \pi^{\perp}x, A \pi^{\perp}x\rangle|\le
16v_3\|\pi^{\perp}x\|_X^2\le 16v_3\|x\|_X^2.$$
It follows that
\begin{equation}\label{4.d}
|F(x)-2\langle \pi x,A\pi^{\perp}x\rangle|\le
48v_3\|x\|_X^2.\end{equation}

We now define $P:=\pi+2\pi A\pi^{\perp}.$ The linear operator $P$ is a
projection onto
$V$; we will use (\ref{4.1}) and so we estimate $\|P\|_Y.$ Assume
$\|y\|_Y=
1$ is chosen so that $\|Py\|_Y=\|P\|.$ Then we may pick $x_1\in V$ with
$d_{Y^*}(x_1,V^{\perp})\le
1$ and
$$ \langle x_1,Py\rangle =\|P\|_Y.$$
Now $x=x_1+\pi^{\perp}(y)\in B_X.$ Note that
$$ F(x)= \langle x_1,\Omega(\pi^{\perp}(y))\rangle= \langle
x_1,\Omega(y)\rangle-\langle x_1,\pi y\rangle.$$
 Hence
$$ |F(x) +\langle x_1,\pi y\rangle| \le 2d_{Y^*}(x_1,V^{\perp})\|y\|_Y
\le 2$$
by Lemma \ref{4.2}.  By (\ref{4.d}) we obtain
$$ |\langle x_1, 2\pi A\pi^{\perp}y +\pi y\rangle|\le 2+48v_3$$
which implies $\|P\|\le 2+48v_3$ and hence gives the estimate $v_3(X)\ge
c\sqrt n$ for suitable $c>0.$\end{proof}

Our second main result of this section gives a rather sharp estimate of
$w_3(\ell_p^n)$ when $p=1$ or $2\le p\le \infty.$ It is a consequence of
more general results which will proved later.

\begin{Thm}\label{4.3}  There are absolute constants $0<c<C<\infty$ so
that for every $n\ge 1$:\newline
(a) $c\log(n+1)\le w_3(\ell_p^n)\le Cp\log(n+1)$ if $p=1$ or $2\le
p<\infty$;\newline
(b) $c\log(n+1)\le w_3(\ell_{\infty}^n)\le C(\log (n+1))^2.$
\end{Thm}

{\it Remark.} We emphasize that $c$ and $C$ are independent of $n$ and
$p.$  We do not have any really good upper estimate for $w_3(\ell_p^n)$
when $1<p<2$, but Theorem 4.3 gives a lower bound in that case:

\begin{Cor}\label{4.4} There is a universal constant $c>0$ so that for
$1<p<2$,
$$ w_3(\ell_p^n)\ge c(p-1)\log (n+1).$$
\end{Cor}

\begin{proof}  We use the following fact proved in an equivalent form in
\cite{21}, p. 21.  There is a universal constant $C$ and for each $n$ a
subspace
$Y_n$ of
$\ell_p^n$,
$1<p<2,$ with
$\dim Y=[n^{2/q}]$ (where $\frac1p+\frac1q=1$) so that:\newline
(a) the Banach-Mazur distance $d(Y_n,\ell_2^{\dim Y_n})\le C;$
\newline
(b) there is a projection $P:\ell_p^n\to Y_n$ with $\|P\|\le C.$

Applying Lemma \ref{4.8} and Lemma \ref{2.3a}  to $Y_n$ we
can find a continuous
2-homogeneous function $f_0:Y_n\to\Bbb R$ with $\omega_3(f_0;B_{Y_n})\le
1$ and $E_3(f_0)\ge c(p-1)\log  (n+1)$  where $c>0$ is a universal
constant.  Defining $f:=f_0\circ P$ we easily have $\omega_3(f)\le C$ but
$E_3(f)\ge c(p-1)\log (n+1)$ and this proves the result.\end{proof}

Except for the case $p=1$, the estimates in Theorem \ref{4.3} will follow
from the following very general estimate:

\begin{Thm}\label{4.5} There are absolute constants $0<c<C<\infty$ so
that for every $n$-dimensional Banach space we have
$$ \frac{c\log (n+1)}{{C_2(X^*)}^{8}}\le w_3(X)\le CT_2(X)^2\log (n+1).$$
\end{Thm}

\begin{proof} {\it (The upper estimate.)} By Theorem \ref{3.11} we have
$w_2(X)\le C(1+\log T_2(X))$ and by Proposition \ref{2.4} we have
$w_3(X)\le C\max(w_2(X),v_3(X)).$  So it will suffice to show a similar
estimate for $v_3(X).$  We obtain the result by a linearization
technique.  We can regard $X$ as $\Bbb R^n$ with an appropriate norm.
Now if $P$ is an $n\times n$ positive-definite matrix, we can define an
$\Bbb R^n-$valued Gaussian random variable $\xi_P$ with covariance
matrix
$P.$ Let $\Gamma$ be the cone of positive-definite matrices.

Suppose now that $f$ is a $2$-homogeneous continuous function on $X$ with
$\omega_3(f;B_X)\le 1.$  We define a function $\hat f$ on $\Gamma$
by putting
$$ \hat f(P):= \mathbf E(f(\xi_P)).$$
Then $\hat f$ is $1$-homogeneous on the cone $\Gamma.$  Let $\Gamma_0$ be
the convex hull of the set of matrices $\{x\otimes x:\ x\in B_X\}$ where
$x\otimes x$ denotes the rank one matrix $(x_ix_j)_{1\le i,j\le n}.$
We need the estimate:

\begin{Lem}\label{4.6} There is a universal constant so that for any
$x_1,x_2\in X$ we have:
\begin{equation}\label{e4.7}
|\frac12(f(x_1+x_2)+f(x_1-x_2))-f(x_1)-f(x_2)|\le C(\|x_1\|^2+\|x_2\|^2).
\end{equation}\end{Lem}

 \begin{proof} By the main result of \cite{2} there is a constant $C_0$
so that $w_3(Y)\le C_0$ for all $2$-dimensional subspaces.  Let $Y:=
\text{ span }\{x_1,x_2\}.$  By Proposition \ref{2.4} there is a quadratic
form $h$ on $Y$ so that $|f(y)-h(y)|\le C\|y\|^2$ for all $y\in Y$ (where
again $C$ is a universal constant).  This immediately yields the
lemma.\end{proof}

Returning to the proof of the theorem we note that if $\xi_P$ and
$\xi_Q$ are independent
then $\xi_P+\xi_Q$ has the same distribution as $\xi_{P+Q}.$  Hence
\begin{equation}\begin{align*}
|\hat f(P+Q)-\hat f(P)-\hat f(Q)| &=|\mathbf E(f(\xi_P+\xi_Q))-\mathbf
E(f(\xi_P)+f(\xi_Q))| \quad\quad\\
= |\mathbf E\frac12(f(\xi_P&+\xi_Q)+f(\xi_P-\xi_Q))-\mathbf
E(f(\xi_P)+f(\xi_Q))|\\
&\le C\mathbf E (\|\xi_P\|^2+\|\xi_Q\|^2).
\end{align*}\end{equation}

Now suppose that $P,Q\in \Gamma_0.$  Then we can write
$P=\sum_{i=1}^ma_ix_i\otimes x_i$ where $\|x_i\|\le 1$ for $1\le i\le m$
and $a_i\ge 0$ with $\sum_{i=1}^ma_i=1.$  Then  $\xi_P$ has the same
distribution as $\sum_{i=1}^ma_i^{1/2}g_ix_i$ where $g_1,\ldots,g_m$ are
independent normalized Gaussian random variables.  Hence
as is well-known (see, e.g., \cite{25}, p.25)
$$\mathbf E(\|\xi_P\|^2)=\mathbf E(\|\sum_{i=1}^ma_i^{1/2}g_ix_i\|^2)\le
T_2(X)^2.$$
Using the similar inequality for $Q$, we obtain
$$ |\hat f(P+Q)-\hat f(P)-\hat f(Q)|\le CT_2(X)^2$$
for a universal constant $C.$  Hence $ \omega_2(\hat f, \Gamma_0)\le
CT_2(X)^2.$  Since $\dim \Gamma_0= \frac12n(n-1)\le n^2 $ we can apply
Theorem \ref{3.1} to $\Gamma_0$ to deduce the existence of an affine
function $\varphi$ on $\Gamma_0$ so that
\begin{equation}\label{e4.8}
|\hat f(P)-\varphi(P)|\le CT_2(X)^2\log (n+1)
\end{equation}
 where $C$ is again a universal constant. In particular
$|\varphi(0)|$ is dominated by $CT_2(X)^2\log(n+1)$ so we can assume that
$\varphi$ is linear on the linear span of $\Gamma_0$.   Let
$h(x)=\varphi(x\otimes x).$  Then $h$ is a quadratic form.  Since $\hat
f(x\otimes x)=\mathbf E(f(gx))=f(x)\mathbf E(g)=f(x)$ where $g$ is a
normalized Gaussian,  we have from (\ref{e4.8})
 $$ |f(x)-h(x)|\le CT_2(X)^2\log (n+1)$$
for all $x\in B_X.$  This gives the desired estimate of $v_3(X)$ and
completes the proof of the upper estimate.

{\it (The lower estimate.)} We establish  a lower
estimate for $v_3(X);$  we first achieve this for the case of
$X=\ell_2^n.$

\begin{Lem}\label{4.8} There is an absolute constant $c>0$ so that for
all $n\ge 1$
\begin{equation}\label{e4.10}
v_3(\ell_2^n)\ge c\log n.
\end{equation}
\end{Lem}

\begin{proof} Let $\varphi(t):=t^2\log |t|$ for $-1\le t\le 1.$  Then,
by the Mean Value Theorem
$$\Delta_h^3\varphi(t) = 3h\Delta_{\theta h}^2\varphi'(t+\theta h)$$
 for some $0<\theta<1.$  Hence $$
|\Delta_h^3\varphi(t)|\le 6\log (1+\sqrt2)|h|^2.$$
Now define for $x\in B_{\ell_2^n}$,
$$ f(x) =\sum_{i=1}^n\varphi(x_i).$$  Then for $x,x+3h\in B_{\ell_2^n},$
$$ |\Delta_h^3f(x)|\le 6\log(1+\sqrt 2)\sum_{i=1}^nh_i^2
<\frac{8}{3}\log(1+\sqrt 2).$$ Hence $
\omega_3(f;B_{\ell_2^n})<
6.$

Since $f$ is even and $f(0)=0$ we can find a quadratic form $h$ on
$\ell_2^n$ so that
$$ \sup_{\|x\|\le 1}|f(x)-h(x)|\le 2E_3(f;B_{\ell_2^n}).$$
As the points $n^{-1/2}\sum_{i=1}^n\epsilon_ie_i\in B_{\ell_2^n}$ for
$\epsilon_i=\pm1$ the left-hand side is at least
$$ \Ave_{\epsilon_i=\pm1}|f
(\frac{1}{\sqrt n}\sum_{i=1}^n\epsilon_i e_i)
-h
(\frac{1}{\sqrt n}\sum_{i=1}^n\epsilon_i e_i)| =|\frac12\log n
+\frac1n\sum_{i=1}^nh(e_i)|.$$
As $f(e_i)=0$ for $1\le i\le n$ we have
$$ \frac1n|\sum_{i=1}^nh(e_i)| \le 2E_3(f:B_{\ell_2}^n).$$
Putting these inequalities together gives
$E_3(f;B_{\ell_2^n})\ge \frac18\log n.$ \end{proof}

Next we need a lemma using the extension constants from Definition
\ref{2.12}.

\begin{Lem}\label{4.7} Let $X$ be an $n$-dimensional Banach space and let
$E$ be a linear subspace of $X.$  Let
${\cal E}_X(E,E^{\perp})=M_1$ and ${\cal E}_X(E,X^*)=M_2.$   Then
\begin{equation}\label{e4.9}
v_3(X/E) \le (M_1+1)(M_2+1)v_3(X).\end{equation}\end{Lem}

\begin{proof} It will be convenient to regard $X$ as $\Bbb R^n$ with an
appropriate norm and let $\langle,\rangle$ be the usual inner-product on
$\Bbb R^n.$ Suppose
$f$ is a
$2$-homogeneous continuous function on
$X/E$
with $\omega_3(f;B_{X/E})\le 1.$  Let $Q:X\to X/E$ be the quotient map.
Then $f\circ Q$ is continuous and $2$-homogeneous on $X$ and
$\omega_3(f\circ Q;B_X)\le 1.$  Hence there is a quadratic form
$h:X\to\Bbb R$ such that
$$ |f(Qx)-h(x)|\le v_3(X)\|x\|_X^2$$
for $x\in \Bbb R^n.$  We can assume $h(x)=\langle x,Ax\rangle $ where $A$
is a symmetric matrix.

Since $\langle x,Ay\rangle=\frac14(h(x+y)-h(x-y))$ we have
$$ |\langle x,Ay\rangle -\frac14(f(Qx+Qy)-f(Qx-Qy))|\le \frac12
v_3(X)(\|x\|_X^2+\|y\|_X^2).$$
Assume $y\in E.$  Then $Qy=0$ and so
$$ |\langle x,Ay\rangle|\le \frac12v_3(X)(\|x\|^2_X +\|y\|_X^2).$$
Replacing $x$ by $\alpha x$ and $y$ by $\alpha^{-1}y$ and minimizing the
right-hand side gives
$$ |\langle x,Ay\rangle|\le v_3(X)\|x\|_X\|y\|_X.$$
This implies that
$$ \|Ay\|_{X^*}\le v_3(X)\|y\|_X$$
when $y\in E.$  From the definition of the extension constant there
exists an $n\times n$ matrix $A_1$ so that $A_1y=Ay$ for $y\in E$ and
$A_1$ has norm at most $M_2v_3(X)$ as an operator from $X$ into $X^*.$
Then $A-A_1$ maps $E$ to $0$ and hence the transpose $A-A_1^t$ maps $\Bbb
R^n$ to $E^{\perp}.$  Now $\|A_1^t\|_{X\to X^*}=\|A_1\|_{X\to X^*}\le
M_2v_3(X)$ and so  $\|A-A_1^t\|_{E\to X^*}\le (M_2+1)v_3(X).$
Using the extension constant again we can find an $n\times n$ matrix
$A_2$ which maps $\Bbb R^n$ into $E^{\perp}$ and such that $\|A_2\|_{X\to
X^*}\le M_1(M_2+1)v_3(X).$

Let $S=A-A_1^t-A_2.$  Then $S$ maps $E$ to $\{0\}$ and $\Bbb R^n$ into
$E^{\perp}.$  It follows that $S=TQ$ where $T$ is a linear operator from
$X/E$ to $E^{\perp}$ and we can define a quadratic form $\psi$ on $X/E$
by $\psi(Qx)=\langle x,Sx\rangle.$

Then \begin{equation}\begin{align*}
 |f(Qx)-\psi(Qx)|&\le |f(Qx)-h(x)|+|\langle x,A_1^tx\rangle|+|\langle
x,A_2x\rangle|\quad\quad\quad\quad \\
&\le (1+ M_2+M_1(M_2+1))v_3(X)\|x\|_X^2\\
&= (M_1+1)(M_2+1)v_3(X)\|x\|_X^2.
\end{align*}
\end{equation}
Now for given $u\in X/E$ we can choose $x\in X$ with $Qx=u$ and
$\|x\|_X=\|u\|_{X/E}.$ This  implies $v_3(X/E)\le (M_1+1)(M_2+1)v_3(X).$
\end{proof}

We can now complete the proof of the lower estimate in Theorem \ref{4.5}.
Suppose $X$ is a Banach space of dimension $n.$  We use the following
powerful form
of the Dvoretzky theorem due to Figiel, Lindenstrauss and Milman \cite{6}
(see \cite {21}, Theorem 9.6, where the theorem is formulated in the form
required here).  There is a subspace $F$ of $X^*$ which is $2$-isomorphic
to $\ell_2^m$ with
\begin{equation}\label{e4.11}
m=\dim F\ge cC_2(X^*)^{-2}n
\end{equation}
We note that the lower estimate in Theorem \ref{4.5} is trivial for
spaces such that
$C_2(X^*)\ge \sqrt{c}n^{1/4}.$  We therefore will consider only those
spaces $X$ for which $C_2(X^*)\le \sqrt{c}n^{1/4}$.  Then
(\ref{e4.11}) gives $m\ge \sqrt n.$

Let us put $E:=F^{\perp}.$  Since $E^*$ is isometric to $X^*/F$ and
$d(F,\ell_2^m)\le 2$ we can apply Theorem 6.9 of \cite{25} to obtain
\begin{equation}\label{e4.12}
C_2(E^*)\le CC_2(X^*)
\end{equation}
where, as usual, $C$ is an absolute constant.

Now we use Corollary \ref{2.14} to estimate the constants $M_1,M_2$ of
Lemma \ref{4.7} as follows:
\begin{equation}\begin{align*}
M_1 &\le
\psi(T_2(X/E))C_2(X^*)(C_2(E^*)C_2(E^{\perp}))^{3/2}\quad\quad\quad\quad\\
M_2 &\le \psi(T_2(X/E))C_2(X^*)(C_2(E^*)C_2(X^*))^{3/2},
\end{align*}\end{equation}
where $\psi:[1,\infty)\to[1,\infty)$ is a suitable increasing function.
Since $X/E$ is isometric to $F^*$ we have $d(X/E,\ell_2^m)\le 2$ and so
$T_2(X/E)\le 2.$  Together with (\ref{e4.12}) this yields
$$ M_1,M_2\le CC_2(X^*)^4.$$
Combining this with (\ref{e4.9})  and Lemma \ref{2.3a} we have:
$$ \frac14
v_3(\ell_2^m) \le v_3(X/E)\le CC_2(X^*)^8v_3(X).$$
Applying now  (\ref{e4.10}) and the inequality $m\ge \sqrt n$ we have
$$ v_3(X)\ge C^{-1}\frac{\log m}{C_2(X^*)^8} \ge
c\frac{\log(n+1)}{C_2(X^*)^8}$$
for an absolute constant $c>0.$
The proof of Theorem \ref{4.6} is now complete.
\end{proof}

{\it Proof of Theorem \ref{4.3}.}  For the case $p=1$ we postpone the
proof to the next section (see Corollary \ref{5.7} below).  For $2\le
p\le \infty$ it suffices to apply Theorem \ref{4.6} to $X=\ell_p^n$
noting that in this case $C_2(X^*)$ is uniformly bounded independent of
$n$ and $p$ while $T_2(\ell_p^n)\le C\sqrt{p}$ for $2\le p<\infty$ and
$T_2(\ell_{\infty}^n)\le C(\log (n+1))^{1/2}$; see, for example
\cite{28}.\qed

\section{Higher order estimates}\label{higher}
\setcounter{equation}0

We now consider upper estimates for $w_n(X)$ when $X$ is a
finite-dimensional Banach space and $n\ge 3$ is arbitrary.  In the proof
we will use heavily the notion and characteristic properties of
$m$-quasilinear functions, which we introduce next:

\begin{Def}\label{5.1} A map $F:X^m\to\Bbb R$ is said to be
{\it $m$-quasilinear} if $F$ is homogeneous in each variable separately
and there is a constant $\lambda\ge 0$ so that for any $1\le j\le m$ and
any $(x_i)_{i\neq j}$ the map
$g_j(x):=F(x_1,\ldots,x_{j-1},x,x_{j+1},\ldots,x_m)$ satisfies
\begin{equation}\label{e5.1}
\omega_2(g_j;B_X) \le \lambda \prod_{i\neq j}\|x_i\|.
\end{equation}
We then set
$\tilde\Delta_m(F)$ to be the infimum of all $\lambda$ so that
(\ref{e5.1})
holds.
\end{Def}

To formulate our main result, we recall that the {\it projection
constant} $\lambda(Y)$ of a finite-dimensional Banach space $Y$ is the
smallest $\lambda\ge 1$ so that if $Y$ is embedded isometrically in a
Banach space $Z$ then there is a linear projection $P:Z\to Y$ with
$\|P\|\le \lambda.$   See for example \cite{Woj}.

\begin{Thm}\label{5.2} For any integers $m\ge k\ge 2$ there is a constant
$C=C(m)$ so that
\begin{equation}\label{e5.2}
w_m(X)\le C\lambda(X^*)^{m-k}w_k(X).
\end{equation}
\end{Thm}

Before proving this estimate we will establish some basic lemmas
on
$m$-quasilinear forms.  We let $C$ denote a constant which depends only
on $m.$

\begin{Lem}\label{5.3}
Suppose $F:X^m\to\Bbb R$ is a symmetric $m$-quasilinear form and that
$f:X\to \Bbb R$ is defined by $f(x)=F(x,\ldots,x).$  Then
\begin{equation}\label{e5.3}
|f(x_1+x_2)-\sum_{k=0}^m\binom{m}{k}F_k(x_1,x_2)|\le
C\tilde\Delta_m(F)\max(\|x_1\|^m,\|x_2\|^m)
\end{equation}
where $F_k(x_1,x_2)=F(x_1,\ldots,x_1,x_2,\ldots,x_2)$ with $x_1$ repeated
$k$ times and $x_2$ repeated $n-k$ times.

More generally there is a constant $C=C(m)$ so that if $x_1,\ldots,x_m\in
X$
\begin{equation}\label{e5.31}
|f(\sum_{i=1}^mx_i)-\sum_{|\alpha|=m}\binom{m}{\alpha}F_{\alpha}(x_1,
\ldots,x_m)|\le C\tilde\Delta_m(F)\max(\|x_1\|^m,\ldots,\|x_m\|^m),
\end{equation}
where we adopt the notation for $\alpha\in \Bbb Z_+^m$ of
$|\alpha|:=\sum_{i=1}^m\alpha_i$ and $$F_{\alpha}(x_1,\ldots,x_m):=
F(x_1,\ldots,x_1,x_2,\ldots,x_2,\ldots,x_m,\ldots,x_m)$$ with each $x_k$
repeated $\alpha_k$ times.
\end{Lem}

\begin{proof} This is established by expanding in each variable
separately and collecting terms.  We omit the details.\end{proof}

Suppose now that $f:X\to\Bbb R$ is a continuous $m$-homogeneous function.
We associate with $f$ the separately homogeneous function $F:X^m\to\Bbb
R$ defined for $\|x_1\|=\|x_2\|=\cdots=\|x_n\|=1$
by
\begin{equation}\label{e5.4}
F(x_1,\ldots,x_m):=\frac{1}{2^mm!}\sum_{\epsilon_i=\pm
1}\epsilon_1\ldots \epsilon_mf(\sum_{i=1}^m\epsilon_ix_i).
\end{equation}
and extended by homogeneity.

\begin{Lem}\label{5.4}
If $f:X\to \Bbb R$ is continuous and $m$-homogeneous then $F$ defined by
(\ref{e5.4}) is symmetric and $m$-quasilinear with $\tilde\Delta_m(F)\le
C\omega_{m+1}(f;B_X).$

Conversely if $F$ is continuous and $m$-quasilinear  then
$f(x):=F(x,\ldots,x)$ is continuous and $m$-homogeneous with
$\omega_{m+1}(f;B_X)\le C\tilde\Delta_m(F).$
\end{Lem}

\begin{proof} Suppose first that $f$ is continuous and $m$-homogeneous
and $F$ is defined by (\ref{e5.4}).  Suppose $(x_i)_{i\neq j}\in B_X$
and $x,x+2h\in B_X.$  Let $E=\spn(\{x_i\}_{i\neq j},x,h).$  Then $\dim
E\le m+1$ and so by the Whitney type result of \cite{2} there is a
constant $C=C(m)$ so that $w_{m+1}(E)\le C.$  By Proposition \ref{2.4} we
also have $v_{m+1}(E)\le C.$  Since $\omega_{m+1}(f;B_E)\le
\omega_{m+1}(f;B_X)$ there is a homogeneous polynomial of degree $m$ on
$E$ so that
$$ |f(u)-g(u)|\le C\|u\|^m\omega_{m+1}(f;B_X)$$
for $u\in E.$  We can express $g$ in the form $g(u)=G(u,\ldots,u)$ where
$G$ is a symmetric $m$-linear form.
Using the polarization formula from multilinear algebra, we have
$$ |F(x_1,\ldots,x_{j-1},u,x_{j+1},\ldots,x_m)-
 G(x_1,\ldots,x_{j-1},u,x_{j+1},\ldots,x_m)| \le C\omega_{m+1}(f:B_X)$$
whenever $\|u\|\le 1$ and $u\in E.$
Let $\phi(u)=F(x_1,\ldots,x_{j-1},u,x_{j+1},\ldots,x_m).$
It now follows that
$$ |\Delta_h^2\phi(x)| \le C\omega_{m+1}(f;B_X)$$
and so
$$ \tilde\Delta_m(F)\le C\omega_{m+1}(f:B_X).$$

We now turn to the converse.  Suppose that $x+jh\in B_X$ for $0\le j\le
m+1.$  Using (\ref{e5.3}) we have
$$ |f(x+jh)-\sum_{k=0}^m\binom{m}{k}j^kF_k(h,x)|\le C\tilde\Delta_m(F).$$
Hence
$$ |\Delta_h^{m+1}f(x)|\le C\tilde\Delta_m(F)$$
as required.
\end{proof}

Our next result shows that symmetric $m$-quasilinear forms can be nicely
approximated by $m$-linear forms.

\begin{Lem}
\label{5.5}
Suppose $F:X^m\to\Bbb R$ is a continuous symmetric $m$-quasilinear form.
Then there is a symmetric $m$-linear form $H:X^m\to\Bbb R$ so that
$$ |F(x_1,\ldots,x_m)-H(x_,\ldots,x_m)|\le
Cv_{m+1}(X)\tilde\Delta_m(F)\prod_{i=1}^m\|x_i\|.$$
\end{Lem}

\begin{proof} Let $f(x):=F(x,\ldots,x).$  By the previous lemma,
$\omega_{m+1}(f:B_X)\le C\tilde\Delta_m(F).$  Hence there is a symmetric
$m$-linear form $H$ so that if $h(x)=H(x,\ldots,x)$ then
\begin{equation}\label{e5.7}
|f(x)-h(x)|\le Cv_{m+1}(X)\tilde\Delta_m(F)\|x\|^m.
\end{equation}
Now let us define $F'$  using (\ref{e5.4}) to be separately homogeneous
and for $\|x_1\|=\|x_2\|=\cdots=\|x_n\|=1,$
\begin{equation}\label{e5.41}
F'(x_1,\ldots,x_n):=\frac{1}{2^mm!}\sum_{\epsilon_i=\pm
1}\epsilon_1\ldots \epsilon_mf(\sum_{i=1}^m\epsilon_ix_i).
\end{equation}
Note that
\begin{equation}\label{e5.6}
\sum_{\epsilon_i=\pm1}\sum_{|\alpha|=m}\binom{m}{\alpha}
\prod_{i=1}^m\epsilon_i^{\alpha_i+1}=m!2^m
\end{equation}
since $\sum_{\epsilon_i=\pm1}\prod_{i=1}^m\epsilon_i^{\alpha_i+1}=0$
unless $\alpha_i=1$ for all $i.$
Hence
$$ \sum_{\epsilon_i=\pm1}\sum_{|\alpha|=m}\binom{m}{\alpha}\epsilon_1
\ldots\epsilon_m F_{\alpha}(\epsilon_1x_1,\ldots,\epsilon_mx_m)=
2^mm! F(x_1,\ldots,x_m).$$
It follows, by Lemma \ref{5.3} that for
$\|x_1\|=\|x_2\|=\cdots=\|x_n\|=1,$
$$ |F'(x_1,\ldots,x_n)-F(x_1,\ldots,x_n)| \le C\tilde\Delta_m(F).$$
We also have, again using Lemma \ref{5.3}
$$ |F'(x_1,\ldots,x_n)-H(x_1,\ldots,x_n)|\le
Cv_{m+1}(X)\tilde\Delta_m(F)$$
and the lemma follows by homogeneity.\end{proof}

{\it Proof of Theorem \ref{5.2}} We will prove by induction that
\begin{equation}\label{ind}
 v_m(X)\le C\lambda(X^*)\max(v_{m-1}(X),v_2(X)) \end{equation}
when $m\ge 3.$

 Let $f:X\to\Bbb R$
be a continuous $m$-homogeneous function with $\omega_{m+1}(f)\le 1.$  We
define $F:X^{m}\to\Bbb R$ using (\ref{e5.4}) so that
$\tilde\Delta_m(F)\le
C.$  Now fixing $u\in X$ we define $$g_u(x):=F(u,x,\ldots,x)$$ so that
$g_u$ is $(m-1)-$homogeneous and $\omega_m(F)\le
C\tilde\Delta_m(F)\|u\|\le
C\|u\|$ by Lemma \ref{5.4}.

Now by Lemma \ref{5.5} there is a symmetric $(m-1)-$linear form
$H_u:X^{m-1}\to\Bbb R$ so that
$$ |F(u,x_2,\ldots,x_m)-H_u(x_2,\ldots,x_m)|\le Cv_{m}(X)\|u\|
\prod_{i=2}^m\|x_i\|.$$
We may clearly suppose that the map $u\to H_u$ is homogeneous.
Then
\begin{equation}\label{e5.8}
 |F(x_1,\ldots,x_m)-H(x_1,\ldots,x_m)|\le
Cv_{m}(X)\prod_{i=1}^m\|x_i\|.\end{equation}

Now let $Z$ be the space of all continuous homogenous functions on $X$
with the norm $\|\varphi\|_Z=\sup_{\|x\|\le 1}|\varphi(x)|.$  Then $X^*$
is a linear subspace of $Z$ and there is a projection $\pi:Z\to X^*$ with
$\|\pi\|\le \lambda(X^*).$

For $x_2,\ldots,x_m\in X$ we define $H_{x_2,\ldots,x_m}$ and
$F_{x_2,\ldots,x_m}\in Z$ by
$$ H_{x_2,\ldots,x_m}(x)=H(x,x_2,\ldots,x_m)$$ and
$$ F_{x_2,\ldots,x_m}(x)=F(x,x_2,\ldots,x_m).$$
Then
$$ d(F_{x_2,\ldots,x_m},X^*)\le
v_2(X)\tilde\Delta_m(F)\prod_{i=2}^m\|x_i\|$$
and by (\ref{e5.8})
$$ \|F_{x_2,\ldots,x_m}-H_{x_2,\ldots,x_m}\|\le
Cv_{m}(X)\prod_{i=2}^m\|x_i\|.$$
Combining we obtain
$$ d(H_{x_2,\ldots,x_m},X^*)\le
C(v_{m}(X)+v_2(X))\prod_{i=2}^m\|x_i\|.$$
Now set
$$ G(x_1,\ldots,x_m)= \pi(h_{x_2,\ldots,x_m})(x_1)$$
so that $G$ is $m$-linear.  Then
$$ |H(x_1,\ldots,x_m)-G(x_1,\ldots,x_m)| \le
(1+\|\pi\|)\prod_{i=1}^m\|x_i\|.$$
Hence appealing again to (\ref{e5.8}) we have
$$ |F(x_1,\ldots,x_m)-G(x_1,\ldots,x_m)|\le
C\lambda(X^*)\max(v_m(X),v_2(X))\prod_{i=1}^m\|x_i\|.$$
This implies (\ref{ind}).

 Since
$v_m(X)\le w_m(X)\le C\max(v_2(X),\ldots,v_m(X)$ by Proposition \ref{2.4}
the theorem is proved.
\qed

\begin{Cor}\label{5.6}For any $m\ge 2$ there  is a constant $C=C(m)$ so
that
$$w_m^{(s)}(n)\le Cn^{\frac{m}{2}-1}\log(n+1)$$ (i.e. for any
$n$-dimensional Banach space  $w_m(X)\le Cn^{\frac{m}{2}-1}\log(n+1).$
\end{Cor}

\begin{proof} Using Theorem  \ref{5.2} with $k=2$ and the
Kadets-Snobar inequality $\lambda(X^*)\le \sqrt n$ (\cite{KS},
\cite{Woj}) we have
$w_m(X)\le
Cn^{\frac{m}{2}-1}w_2(X)$, but $w_2(X)\le C\log(n+1)$ by Theorem
\ref{3.1}.\end{proof}

\begin{Cor}\label{5.7} For any $m\in\Bbb N$ there exists a constant
$C=C(m)$ so that $$C^{-1}\log(n+1)\le w_m(\ell_1^n)\le
C\log(n+1).$$\end{Cor}

\begin{proof} Since $\lambda(\ell_{\infty}^n)=1$ (see e.g \cite{26}) by
Theorem
\ref{5.2} with
$k=2$ we have $w_m(\ell_1^n)\le Cw_2(\ell_1^n)\le C\log(n+1).$
Conversely by Corollary \ref{2.5} and Theorem \ref{3.8} we have
$C^{-1}\log(n+1)\le w_2(\ell_1^n)\le Cw_m(\ell_1^n).$\end{proof}

\begin{Cor}\label{5.8}  For any $m\ge 3$ and $2\le p<\infty$ there is a
constant $C=C(m,p)$ so that
$$ w_m(\ell_p^n)\le Cn^{\frac{m-3}{2}}\log(n+1).$$
\end{Cor}

\begin{proof} Apply Theorem \ref{5.2} with $k=3$ and use Theorem
\ref{4.2}.\end{proof}

There is a striking difference between the results for $p\ge 1$ and for
$0<p<1$, when the sets $B_{\ell_p^n}$ are no longer convex.  The
following Theorem is then true:

\begin{Thm}\label{5.9}
If $0<p<1$ and $m\ge 2$ there is a constant $C=C(p,m)$ so that
$w_m(\ell_p^n)\le C$ for all $n\ge 1.$\end{Thm}

\begin{proof}  It is easily checked that the proof of Theorem \ref{5.2}
goes through with trivial changes for $r$-normed spaces when $r<1$ (see
Remark after Corollary \ref{2.5}).  Of course the constant $C$ in its
formulation depends
now on $r.$  Applying this result to $\ell_p^n$ with $r=p<1$ we therefore
have
$$w_m(\ell_p^n)\le C(m,p)\lambda((\ell_p^n)^*)^{m-2}w_2(\ell_p^n).$$
But $(\ell_p^n)^*=\ell_{\infty}^n$ and it is essentially proved in
\cite{9} (in an equivalent formulation related to the notion of a $\cal
K$-space) that $w_2(\ell_p^n)\le C(1-p)^{-1}$ with $C$ an absolute
constant independent of $n.$  This proves the Theorem.\end{proof}

\end{document}